\documentclass[12pt]{amsart}
\usepackage{multirow}
\usepackage[margin=3cm]{geometry}
\usepackage{amsmath, amsthm, amsfonts}
\usepackage{verbatim}
\usepackage{natbib, float, subcaption}
\usepackage{bbm}
\usepackage[colorlinks=TRUE,citecolor=blue]{hyperref}
\usepackage{algorithm}
\usepackage{multirow}
\usepackage{bm}
\newcommand*{\B}[1]{\ifmmode\bm{#1}\else\textbf{#1}\fi}
\usepackage{pdfpages}

\usepackage{xcolor}
\usepackage{enumitem}
\usepackage{datetime}
\usepackage{mathtools}
\usepackage{graphicx}
\graphicspath{ {./images/} }
 \usepackage{algorithm}
\usepackage{algpseudocode}
\usepackage{comment}
\usepackage{placeins}
\newtheorem{theorem}{Theorem}
\newtheorem{lemma}{Lemma}

\newcommand{\E}{\mathbb{E}}
\newcommand{\var}{{\rm var}}

\title[Bayesian RDD with heterogeneous treatment effects]{Bayesian analysis of regression discontinuity designs with heterogeneous treatment effects \\ Last update: \shortdate\today \ at \currenttime \bigskip }

\author{Kevin Tao, Y. Samuel Wang, David Ruppert}

\begin{document}
\begin{abstract}
Regression Discontinuity Design (RDD) is a popular framework for estimating a causal effect in settings where treatment is assigned if an observed covariate exceeds a fixed threshold. We consider estimation and inference in the common setting where the sample consists of multiple known sub-populations with potentially heterogeneous treatment effects. In the applied literature, it is common to account for heterogeneity by either fitting a parametric model or considering each sub-population separately. In contrast, we develop a Bayesian hierarchical model using Gaussian process regression which allows for non-parametric regression while borrowing information across sub-populations. We derive the posterior distribution, prove posterior consistency, and develop a Metropolis-Hastings within Gibbs sampling algorithm. In extensive simulations, we show that the proposed procedure outperforms existing methods in both estimation and inferential tasks. 
Finally, we apply our procedure to U.S. Senate election data and discover an incumbent party advantage which is heterogeneous over different time periods.
\end{abstract}

\maketitle


\section{Introduction}
The regression discontinuity design (RDD) was proposed by \cite{original} as a quasi-experimental framework to estimate a causal effect 
in settings where treatment assignment depends only on whether a specific covariate---often referred to as the running variable---exceeds a fixed threshold or cut-off. 
Under certain assumptions, the local average treatment effect---i.e., the causal effect for individuals with a running variable at the cut-off---can be identified and estimated by, roughly speaking, comparing the conditional expectation of the outcome 
just above and just below the cut-off~\citep{imbens2008regression,Guide}. 
Regression discontinuity designs are commonly used in various disciplines such as political science, public policy, economics, and medicine. For instance, a county-level income cut-off in eligibility for government funding is used to estimate the causal effect of the ``Head Start'' program on child mortality~\citep{child_mortality}; an age based cut-off in Medicare eligibility is used to estimate the effect of health insurance on health provider utilization habits \citep{medicare}; and an income cut-off in financial aid eligibility is used to estimate the effect of financial aid on college enrollment \citep{Finaid}. 

In this paper, we focus on RDD settings where the population comprises known sub-populations with potentially heterogeneous local average treatment effects. Our goal is to estimate and conduct inference on the local average treatment effect for each sub-population. 
This is a common goal in the applied literature which is typically achieved by either adding an interaction term or fitting separate non-parametric regressions for each sub-population. However, because each sub-population may have a small sample size, these approaches typically suffer from imprecise estimation and confidence intervals which do not cover at the nominal rate.

As an alternative, we propose a Bayesian hierarchical modeling approach using Gaussian process regression. 
We establish posterior consistency for our estimation procedure; and most notably, the proposed procedure shows very strong empirical results when compared to previously proposed procedures. We provide extensive simulations which show that the point estimator has smaller mean squared error than existing methods. In addition, the credible intervals attain frequentist coverage properties and are much shorter than existing methods.  



The remainder of this section provides an overview of the standard RDD setup and existing literature. Section 2 introduces the setup for RDD with sub-populations, along with our model and its posterior derivation and posterior consistency theorem. Section~3 presents a Metropolis within Gibbs sampler for our model. Section 4 contains numerical experiments, and Section 5 contains an analysis of incumbency advantage in U.S. Senate elections. 

\subsection{Background and literature review}\label{sec:background}
In the standard RDD setting, we observe $(T_i,Y_i,Z_i)$ for the $i$th unit where $T_i$ is a treatment indicator which is $1$ if the $i$th unit receives treatment and is $0$ otherwise, $Y_i$ is the observed outcome, and $Z_i$ is the running variable. 
Using the potential outcome framework~\citep{Rubin1974EstimatingCE, holland}, let $Y_i(0)$ denote the outcome we would have observed for the $i$th observational unit if it had not received treatment and $Y_i(1)$ denote the outcome we would have observed if it had received treatment. The observed outcome can be written as $Y_i = T_iY_i(1)+(1-T_i)Y_i(0)$.

The key assumption for RDD is that the conditional probability of receiving treatment, $P(T_i=1\mid Z_i=z)$, takes a discontinuous jump at the cut-off threshold, which we can assume without loss of generality to be $z = 0$. 
The fuzzy RDD setting only requires that $\lim_{z\rightarrow 0^-} P(T_i=1|Z_i=z) \neq \lim_{z\rightarrow 0^+} P(T_i=1|Z_i=z)$. In contrast, the sharp RDD setting---which we consider in this paper---requires that every individual above the cut-off receives treatment and no one below the cut-off receives treatment; i.e., $P(T_i \mid Z_i < 0) = 0$ and $P(T_i \mid Z_i\geq 0) = 1$.  
Because the running variable may also effect the outcome, the average treatment effect, $\E\left(Y_i(1) - Y_i(0)\right)$, is not identifiable in a sharp RDD due to the lack of overlap in the running variable between the treated and control groups. Nonetheless, when the conditional expectation of the potential outcomes is continuous with respect to $Z$, the local average treatment effect (LATE)---i.e., $    \delta = \E(Y_i(1)-Y_i(0)|Z_i=0)$ which is the treatment effect for units at the threshold---can be identified~\citep{identify}.  




There are two common frameworks for estimation and inference: the local randomization approach and the the local continuity approach~\citep{Guide}. The local randomization approach assumes that within a neighborhood around the threshold, the running variable---and thus the treatment---is assigned ``as if random'' so that the conditional expectations of the potential outcomes are constant within the neighborhood~\citep{senate}. 
In contrast, we use the local continuity approach which only requires the conditional expectation of the potential outcomes to be continuous within a neighborhood on each side of the threshold. The local ATE can then be identified by 
\begin{equation}
\begin{aligned}\label{eq:lr_limits}
    \delta &= \E(Y_i(1) \mid  Z_{i}=z) - \E(Y_i(0)\mid Z_{i}=z)\\
    &= \lim_{z\rightarrow 0^+}\E(Y_i \mid Z_{i}=z) - \lim_{z\rightarrow 0^-}\E(Y_i \mid  Z_{i}=z).
\end{aligned}
\end{equation}

The right and left limits in Eq.~\eqref{eq:lr_limits} are typically estimated by applying local linear regression (LLR)~\citep{Fan1994LocalPM} on both sides of the cut-off. 
However, the performance of LLR depends heavily on the bandwidth selection. \citet{IK} propose a data-driven procedure for selecting the mean squared error optimal bandwidth. \citet{RBC} develop an alternative robust bias-corrected method which yields improved confidence interval coverage, and~\citet{ce_opt} propose a procedure for selecting the ``coverage optimal'' bandwidth. As Bayesian alternatives to LLR,  \citet{chib_greenberg_simoni_2022} use Bayesian cubic splines and \citet{GPR_branson} use Gaussian process regression. Similar to the frequentist approach, both of these methods fit separate regressions for the treated and control groups, and then estimate the treatment effect by extrapolating to the cut-off. Similar to~\citet{GPR_branson}, we also use Gaussian process regression to fit the conditional expectations. However, our approach differs in two significant ways. Most notably, we specify a hierarchical model similar to multi-task Gaussian process regression \citep{Leroy_2022} which allows for borrowed information across sub-populations. In addition, we use a different regression specification which fits the conditional expectation on both sides of the cut-off simultaneously. 


Investigation of heterogeneous treatment effects is common in the applied literature; e.g., \cite{obese} find that the effect of body weight on labor income in Mexico varies across ethnic group and gender and \cite{RD_invest} assess how government fund affects firms' willingness to conduct R\&D, and found the effect to vary across the amount of capital originally owned by the firm. However, the methodological literature for estimating heterogeneous treatment effects for RDD is less developed. \citet{becker} account for heterogeneity by applying LLR with an added interaction term. 
\cite{HSU2019468} consider a setting where the treatment effect may vary with respect to other covariates and propose a global test for effect heterogeneity. 
\cite{reguly2021heterogeneous} considers the problem of discovering sub-groups over which heterogeneity occurs. They propose using a regression tree to first partition the covariate space, and subsequently estimate a treatment effect within each leaf using polynomial regressions. For valid inference, this requires splitting the data into three sets for training, testing, and estimation. \citet{alcantara2024modified} propose using a modified Bayesian additive regression tree for the same purpose. 


Most similar to our work, \cite{sugasawa2023hierarchical} also consider the setting where the sub-populations are known a priori. They propose a Bayesian hierarchical pseudo-likelihood model using LLR to estimate the treatment effect. Our procedure also uses a hierarchical model, but instead employs Gaussian process regression. This comes at a computational cost, but---as shown in the numerical results---yields substantially smaller mean squared error and shorter credible intervals.

\section{Model and methodology}
We generalize the standard RDD setting described in Section~\ref{sec:background}, by allowing the population to be composed of $J$ known sub-populations. In the $j$th sub-population, we have $n_j$ units, and we observe $(Y_{ij},Z_{ij},T_{ij})$ for each unit $i = 1,\ldots,n_j$. As before, $T_{ij}$ is the treatment indicator,  $Z_{ij}$ is the running variable, and $Y_{ij}$ is the observed outcome. Additionally, $Y_{ij}(0),Y_{ij}(1)$ denote the potential outcomes so that $Y_{ij} = T_{ij}Y_{ij}(1) + (1-T_{ij})Y_{ij}(0)$. Most notably, we allow the conditional expectations $\E(Y_{ij}(1)\mid Z_{ij})$ and $\E(Y_{ij}(0)\mid Z_{ij})$ to differ across sub-populations. Thus, the estimand of interest, the local ATE for the $j$th sub-population, denoted as $\delta_j = \E(Y_{ij}(1) - Y_{ij}(0) \mid Z_{ij} = 0)$, may also vary across sub-populations.
As expected, when each sub-population satisfies the conditions in~\citet{identify}, each sub-population LATE is identifiable. We formalize this in Lemma~\ref{lem:id} below. 

\begin{lemma}\label{lem:id}
Suppose for all $j=1,\dots,J$:
    \begin{enumerate}
    \item $T_j^+ := \lim_{z \rightarrow 0^+} \E(T_{ij}|Z_{ij}=z) \neq \lim_{z \rightarrow 0^-} \E(T_{ij}|Z_{ij}=z) =: T_j^-$
    \item $\E(Y_{ij}(0) |Z_{ij}=z)$ is continuous at $z=0$
    for all $j=1,\dots,J$
    \item $\E(Y_{ij}(1)|Z_{ij}=z)$ is continuous at $z=0$ for all $j=1,\dots,J$
    \item $T_{ij}\perp\!\!\! \perp \delta_{ij} | Z_{i,j}$,
\end{enumerate}
with $\delta_{ij} = Y_{ij}(1)-Y_{ij}(0)$. Then the treatment effect for each sub-population is identifiable.  
\end{lemma}
In a sharp RDD, $T_j^+=1$ and $T_j^- =0$, so (1) holds automatically.  Although we will assume a sharp RDD in the remainder of the paper, in Lemma~\ref{lem:id}  we allow a fuzzy RDD so assumption (1) is necessary.

\subsection{Model}
Let ${\rm GP}(m, K)$ denote a Gaussian process with mean function $m$ and covariance kernel $K$, and let $\mathbf{0}(z = \cdot)$ denote the function over $z$ which is $0$ everywhere. 
The estimand of interest is the local ATE for each sub-population, $\delta_j = \E(Y_{ij}(1) \mid Z_{ij} = 0) - (\E(Y_{ij}(0) \mid Z_{ij} = 0)$. 
Furthermore, let
\[f_j(z) = 1_{\{z < 0 \}}\E(Y_{ij}(0) \mid Z_{ij} = z) + 1_{\{z\geq 0\}}[\E(Y_{ij}(1) \mid Z_{ij} = z) - \delta_j].   \]
We assume that outcomes are drawn from the following generative model: 
\begin{equation}
\begin{aligned}\label{eq:HGPR}
        Z_{ij} &\sim P_{Z} \\
        g(z=\cdot)&\sim {\rm GP}\left(\mathbf{0}(z = \cdot),K_g\right), \\
        f_j(z=\cdot)&\sim {\rm GP}\left(g,K_j\right) \quad  \text{ for } j = 1, \ldots, J, \\
         (\delta_1, \ldots, \delta_J) &\sim N\left(\mu,K_{\delta}\right),\\
         \epsilon_{ij} & \sim N\left(0, T_{ij}\sigma^2_{+j} + (1- T_{ij})\sigma^2_{-j}\right), \\
        Y_{ij} &= f_j(Z_{0ij}) + T_{ij} \delta_j +\epsilon_{ij} \quad \text{ for }  i = 1, \ldots, n_{j}.
\end{aligned}
\end{equation}
We also suppose that each $\sigma_{+j}$ and $\sigma_{-j}$ are drawn i.i.d. from priors $\nu_+$ and $\nu_+$ respectively.

In the model above, we assume that all $f_j$ have a common prior mean $g$, which itself is a mean $0$ Gaussian process. This is similar to the multi-task Gaussian process formulated in~\citet{Leroy_2022} and enforces similarity across the general shape of each $f_j$. The mean $0$ prior on $g$ could be generalized to a polynomial with coefficients drawn from some prior distribution. One straightforward choice would be polynomial with degree 1, which---as discussed in \citet{GPR_branson}---is the Bayesian analog to LLR. However, \citet{GPR_branson} found that the empirical results were insensitive to the choice of mean function; thus, for simplicity, we will proceed with the mean $0$ model. 

By definition, $f_j$ is continuous across the threshold; furthermore, in the generative model, $f_j$ is a Gaussian process, which implies that the derivatives of $\E(Y_{ij}(0) \mid Z= z)$ and $\E(Y_{ij}(1) \mid Z= z)$---if they exist---are equal at $z = 0$. 
However, we emphasize, that this does not constrain the treatment effect to be constant with respect to $Z_{ij}$. Moreover, the numerical experiments show that our procedure outperforms existing procedures even when the assumption of equal derivatives at 0 does not hold.

We also assume that the treatment effects $(\delta_1, \ldots, \delta_J)$ are drawn from a Gaussian distribution with covariance $K_\delta$. In most settings where the analyst does not have strong prior knowledge about which groups may have similar treatment effects, $K_\delta$ should be set so that the $\delta_j$ are independent.
However, in settings where the similarity between groups may be quantified or some groups are expected to be more similar than others,  the model can  incorporate this prior information flexibly.   
In Section~\ref{sec:posterior}, we derive the posterior for the treatment effects in the most general case. However, Theorem~\ref{thm:hgpr_consistency} focuses on the default setting where the $\delta_j$ are independent in order to analyze posterior consistency.

Given the model in Eq.~\eqref{eq:HGPR}, the parameters of interest are \[(f_{j}, \sigma^2_{+j},\sigma^2_{-j},\delta_j  \,:\, j = 1, \ldots, J)\in \mathcal{H}^J \times (\mathbb{R}^+)^{2J}\times \mathbb{R}^J,\]
where $\mathcal{H}^J$ is the Cartesian product of function spaces for $(f_1, \ldots, f_J)$. Theorem~\ref{thm:hgpr_consistency} describes conditions under which the posterior concentrates around the true parameters; it assumes that $g$ and the kernels $K_j$ are known and that each $\delta_j$ is drawn i.i.d. from some distribution $P_\delta$.  We will use $(f_{j0}, \sigma^2_{+j0},\sigma^2_{-j0},\delta_{j0} )$ to denote the true parameters, and let $\Pi$ and $\Pi(\cdot \mid \mathbf{Y},  \mathbf{Z})$ denote the joint prior and posterior respectively. We also use the metric $d_{P_Z}(f_1,f_2)=\inf \{\epsilon:P_Z(\{x:|f_1(x)-f_2(x)|>\epsilon\})<\epsilon\}$
, with ${P_Z}$ being the distribution of the running variable.

\begin{theorem}[Posterior Consistency] \label{thm:hgpr_consistency}
Suppose the following assumptions hold:
\begin{enumerate}
    \item The running variable $Z_{ij}\stackrel{iid}{\sim} P_Z$, where $P_Z$ is defined on $[-1,1]$ and continuous.
    \item\label{asm:kernel} The Gaussian processes $f_j$ have a continuously differentiable mean function $g$ and the kernel function $K_j(z,z')$ has continuous 4th partial derivatives. The function $g$ and kernel $K_j$ are known.
    \item\label{asm:prior} The priors for $\sigma_{j+},\sigma_{j-}, \delta_{j}$---which we denote as $\nu_+,\nu_-,P_{\delta}$ respectively---assign positive probability to every neighborhood of the true values $\sigma_{+j0},\sigma_{-j0},\delta_{j0}$. Furthermore, the distribution $P_{\delta}$ is sub-exponential; i.e. there exist $K>0$ such that $P(|\delta_j|\geq t) = O(e^{-Kt})$.
\end{enumerate}
Let {\footnotesize \begin{equation*} 
    U_{j\epsilon} = \left\{(f_j,\sigma_{+j},\sigma_{-j},\delta_j):d_{P_Z}(f_j,f_{j0})<\epsilon,|\sigma_{+j}/\sigma_{+j0}-1|<\epsilon,|\sigma_{-j}/\sigma_{-j0}-1|<\epsilon,|\delta_j-\delta_{j0}|<\epsilon\right\}.
\end{equation*}
}
$U_{\epsilon} = \bigtimes\limits_{j=1}^J U_{j\epsilon}$, where $A\bigtimes B$ denotes the Cartesian product of the sets $A,B$. Then we have that \begin{equation*}
    \Pi(U^C_{\epsilon}|\B{Y},\B{Z})\xrightarrow{}0, \quad a.s.\; [P_{\theta_0}].
\end{equation*}
\end{theorem}
The proof of Theorem~\ref{thm:hgpr_consistency} is contained in the appendix. Assumption~\ref{asm:kernel} on the kernel holds for the square exponential kernel, and the continuous differentiability of the function $g$ is automatic when $g(z=\cdot)\sim {\rm GP}\left(\mathbf{0}(z = \cdot),K_g\right)$ with $K_g$ being again the square exponential kernel.
The normality assumption for $\delta_j$ satisfies the sub-exponential assumption. Furthermore, commonly used priors, such as Gamma or half-Cauchy, also satisfy Assumption~\ref{asm:prior} on $\nu_-,\nu_+$. 

\subsection{Posterior distribution of $\B{\delta}$}
\label{sec:posterior}
We now derive the posterior distribution of $\B{\delta}=[\delta_1,\dots,\delta_J]$.  We first introduce some additional notation.
\begin{itemize}
    \item Let $\mathbf{Y}_{+j} = (Y_{ij} \, : T_i = 1)$, $\mathbf{Y}_{-j} = (Y_{ij} \, : T_i = 0)$,  $\mathbf{Y}_- = [\mathbf{Y}_{-1}, \mathbf{Y}_{-2}, \dots, \mathbf{Y}_{-J}]$, $\mathbf{Y}_+ = [\mathbf{Y}_{+1}, \mathbf{Y}_{+2}, \dots, \mathbf{Y}_{+J}]$, and  $\mathbf{Y} = [\mathbf{Y}_-,\mathbf{Y}_+]$.
    Let $n_{+j}$ and $n_{-j}$ be the dimensions of $\mathbf{Y}_{+j}$ and $\mathbf{Y}_{-j}$ respectively so that $n_j = n_{+j} + n_{-j}$. Let $N = \sum_j^J n_j$, $N_+ = \sum_j^J n_{+j}$, and $N_- = \sum_j^J n_{-j}$.  
    For the running variable, we let $\mathbf{Z}_{+j}$, $\mathbf{Z}_{-j}$, $\mathbf{Z}_{+}$, $\mathbf{Z}_{-}$, and $\mathbf{Z}$  be defined analogously.
    \item 
    Let $\B{g} = [g(\mathbf{Z}_-),g(\mathbf{Z}_+)]$, and let $\B{K_g}$ be the kernel matrix obtained by evaluating the covariance kernel $K_g$ on $\mathbf{Z}$. Similarly, let $\mathbf{f}_-$ and $\mathbf{f}_+$ denote the vector formed by evaluating the appropriate $f_j$ on $\mathbf{Z}_-$ and $\mathbf{Z}_+$ respectively, and let $\mathbf{f} = [\mathbf{f}_-, \mathbf{f}_+]$.  
    \item 
    $\mathbf{H} \in \{0,1\}^{N_+\times J}$, where the $\mathbf{H}_{i, j}$ is $1$ if the $i$th element of $\mathbf{Z}_+$ is from the $j$th sub-population and $\mathbf{H}_{i, j} = 0$ otherwise.
    \item $\B{K_{\delta}}$ is the $J \times J$ matrix obtained by applying $K_{\delta}$ to $(1, \ldots, J)$.
    \item $\mathbf{D}$ is a $N\times N$ matrix. If the $k$th and $l$th element of $\mathbf{Z}$ both belong to the same sub-population, $j$, then $\mathbf{D}_{k,l} = K_j(\mathbf{Z}_{k},\mathbf{Z}_{l})$; otherwise $\mathbf{D}_{k,l} = 0$.
    \item $\mathbf{\Sigma}$ is a block diagonal matrix where the blocks are $\sigma_{-1}^2 I_{n_{-1}}, \ldots,  \sigma_{-J}^2 I_{n_{-J}}$ and  $\sigma_{+1}^2 I_{n_{+1}}, \ldots,  \sigma_{+J}^2 I_{n_{+J}}$, where $I_n$ denote the identity matrix of dimension $n$  
\end{itemize}
We also state the following property of the multivariate Gaussian distribution which will be useful for our derivation; results of this form can be found in \cite{cond_normal}. If $X\sim N(\mu,\Sigma)$, and $Y|X \sim N(AX+b,\Lambda)$, then:

\begin{equation}\label{useful}
\begin{pmatrix}
    X\\
    Y
\end{pmatrix} \sim N\left( \begin{pmatrix}
    \mu \\
    A\mu+b
\end{pmatrix},\begin{bmatrix} 
\Sigma & \Sigma A^T  \\
A\Sigma & A\Sigma A^T+\Lambda 
\end{bmatrix}\right).
\end{equation}

Throughout this section, we will condition on $\mathbf{Z}$, though for notational simplicity we will not state the conditioning explicitly. For now, we will also assume that the error variances as well as parameters in all the kernel functions $K_1,\dots,K_J,K_g,K_{\delta}$ are fixed hyperparameters.

First we observe that $\mathbf{f}|\B{g}\sim N(\B{g},\mathbf{D})$, and marginalizing out $\B{g}$ results in $\mathbf{f} \sim N(\B{0},\B{K_g}+\mathbf{D})$. Furthermore, because $\B{\delta} \perp\!\!\!\perp \mathbf{f}$, we have:
\begin{equation}
    \begin{pmatrix}
    \mathbf{f} \\
    \bm{\delta}
\end{pmatrix}\sim N\left( \begin{pmatrix}
    \B{0} \\
    \B{1}\mu
\end{pmatrix},\begin{bmatrix} 
\bm{K_{g}}+\mathbf{D} & \bm{0}  \\
\bm{0} & \bm{K_{\delta}} 
\end{bmatrix}\right),
\end{equation}
and conditional on $(\B{\delta},\mathbf{f})$, we have:
\begin{equation}
    \B{Y}| \mathbf{f},\bm{\delta}, \sim N\left( \begin{pmatrix}
    \mathbf{f}_0 \\
    \mathbf{f}_1+\mathbf{H}\B{\delta}
\end{pmatrix},\mathbf{\Sigma} \right).
\end{equation}
Marginalizing out $\mathbf{f}$, results in the joint distribution:
\begin{equation}\label{joint}
    \begin{pmatrix}
    \bm{\delta} \\
    \B{Y}
\end{pmatrix}\sim N\left( \begin{bmatrix}
    \mu_{\bm{\delta}} \\
    \mu_{\bm{Y}}
\end{bmatrix},\Lambda \right),
\end{equation}
with
\begin{equation*}
    \begin{split}
        \mu_{\bm{\delta}} &= \bm{1}\mu \\
        \mu_{\bm{Y}} &= \begin{pmatrix}
    \B{0} \\
    \B{1}\mu
    \end{pmatrix} \\
    \Lambda &= \begin{bmatrix} 
\bm{K_{\delta}} & \begin{pmatrix}\B{0} & \bm{K_{\delta}}\mathbf{H}^T \end{pmatrix}  \\
\begin{pmatrix} \B{0} \\ \mathbf{H}\bm{K_{\delta}} \end{pmatrix} & 
\bm{K_{g}}+\mathbf{D}+\mathbf{\Sigma}+ 
\begin{bmatrix} 
\B{0} & \B{0}  \\
\B{0} & \mathbf{H}\bm{K_{\delta}}\mathbf{H}^T
\end{bmatrix}
\end{bmatrix}.
    \end{split}
\end{equation*}
This results in a fully analytic posterior:
\begin{equation}
    [\B{\delta}|\B{Y}]\sim N(\bm{1}\mu+\Lambda_{12}\Lambda_{22}^{-1}(\bm{Y}-\mu_{\bm{Y}}),\Lambda_{11}-\Lambda_{12}\Lambda_{22}^{-1}\Lambda_{21})
\end{equation}
when the hyperparameters are known, and $\Lambda_{ij}$ denotes the block matrices in $\Lambda$

\section{Computation}
For each $f_j$, we use the squared exponential covariance kernel: $K_j(x,y) = r_j^2\exp(-\frac{1}{2l_{j}^2}(x-y)^2)$. We will further assume that the same kernel parameters are used across the mean functions of all sub-populations; i.e., $r_1 = \ldots = r_J$ and $l_1 = \ldots = l_J$.
Note that this would be similar to using the same bandwidth for each sub-population in a LLR analysis. With this simplification there are a total of 6 hyperparameters in the kernel functions $K_1,K_g,K_{\delta}$ and $2J$ error variances, in addition to $\mu$. We will let $\theta = \left\{\mu,\{\sigma^2_{+j}\}_{j=1}^J,\{\sigma^2_{-j}\}_{j=1}^J,r_{\delta}^2,r_1^2,r_g^2,1/l_{\delta}^2,1/l_1^2,1/l_g^2\right\}$. 

\subsection{Sampling from the posterior}

We will take a fully Bayesian approach and place an independent half Cauchy prior on each element of $\theta$ except $\mu$, which is handled with a diffuse normal prior. 
The half Cauchy prior does not yield an analytic conditional distribution and the Metropolis-Hastings acceptance rate may be low in the $2J+7$-dimensional parameter space. Thus, we propose a Metropolis-Hastings within Gibbs sampling method to iterate over every element of $\theta$ using Metropolis updates through the joint density $[\B{Y},\theta]$. Following the practice in Bayesian statistics literature, we will use the following notations for densities of random variables:
\begin{equation*}
    \begin{split}
        &[X] = \text{density of $X$}\\
        &[X,Y] = \text{joint density of $X,Y$}\\
        &[X|Z] =\text{conditional density of $X|Z$}.
    \end{split}
\end{equation*}  We will use $m_j(\cdot|\cdot)$ to denote the proposal density for the $j$th element of $\theta$. Since elements of $\theta_{2:2J+7}$ are restricted to $(0,\infty)$, we can apply a normal proposal to the log scale, effectively leading to a lognormal proposal $\theta_j'\sim \mbox{lognormal}(\log(\theta_j),\sigma^2)$. $\theta_1$ will be given a normal prior and normal proposal. The details are given in Algorithm~\ref{hgpr_mcmc}.

\begin{algorithm}
\caption{HGPR MCMC}\label{hgpr_mcmc}
\begin{algorithmic}
\Require $\B{Y},\B{Z},T,\{m_j(\cdot|\cdot)\},[\theta]$
\State draw $\theta^0$ from the prior $[\theta]$
\State $t \gets 1$
\While{$t \leq T$}
\For{$j=1,\dots,2J+7$}
    \State Generate a new $\theta_j'$ from $m_j(\theta_j'|\theta_j^{t-1})$
    \State Calculate acceptance probability $a = \min \left(\frac{[\B{Y},\theta_{1:j-1}^t,\theta_j',\theta_{j+1:2J+7}^{t-1}]m_j(\theta_j^{t-1}|\theta_j')}{[\B{Y},\theta^{t-1}]m_j(\theta_j'|\theta_j^{t-1})},1 \right)$
    \State With probability $a$ set $\theta_j^t = \theta_j'$, else $\theta_j^t = \theta_j^{t-1 }$ 
\EndFor
\State Draw $\B{\delta}$ from $[\B{\delta}|\theta^t,\B{Y}]$
\State $t\gets t+1$
\EndWhile
\State Output $\{\theta^t,\delta_1^t,\dots,\delta_J^t\}_{t=1}^T$
\end{algorithmic}
\end{algorithm}

After drawing a sufficient number of samples from the posterior, we can use the MCMC sample mean as point estimates, and MCMC empirical posterior intervals for inference. At each iteration $t$, we can further draw $\B{\delta}$ from its full conditional $[\B{\delta}|\theta^{t},\B{Y}]$. 
This yields posterior draws from the  conditional distribution of distribution $\delta,\theta \mid \mathbf{Y}$. 

\subsection{MCMC output analysis}
We use the posterior mean $\frac{1}{T}\sum_{t=1}^T\B{\delta}^t$ to estimate the treatment vector $\bm{\delta}$, and use the MCMC sample quantiles to create marginal posterior credible intervals for component-wise inference.

To create the simultaneous credible region, we first note that the posterior $\B{\delta}|Y$ is not analytic due to randomness in $\theta$. Since components of $\bm{\delta}^t$ are not independent, we can continue to exploit correlations in $\B{\delta}$ to create small confidence region. To do so, we define the critical value $R_{\alpha}$ as:
\begin{equation}\label{post_elip}
        R_\alpha = \min\left\{\, R \;: \; \frac{1}{T}\sum_{t=1}^T\B{1}\left\{(\B{\delta}^t-\hat{\mu}_{\B{\delta}})^T\hat{\Sigma}_{\B{\delta}}^{-1}(\B{\delta}^t-\hat{\mu}_{\B{\delta}})<R\right\}\geq 1-\alpha \right\},  
\end{equation}
where $\hat{\mu}_{\B{\delta}}$ is the posterior mean vector, and $\hat{\Sigma}_{\B{\delta}}$ is the batch mean estimator from \citep{vats2019multivariate}, which is known to be a strongly consistent estimator for the posterior covariance matrix in MCMC. Finally, we can define our confidence region as:
\begin{equation}\label{post_elip}
    C(\alpha) = \{\delta:(\delta-\hat{\mu}_{\B{\delta}})^T\hat{\Sigma}_{\B{\delta}}^{-1}(\delta-\hat{\mu}_{\B{\delta}})<R_{\alpha}\},
\end{equation}
which has a volume given by the formula:
\begin{equation}
    {\rm Vol}(C(\alpha)) = \frac{2\pi^{p/2}}{p\Gamma(p/2)}R^{p/2}|\hat{\Sigma}_{\B{\delta}}|^{1/2}.
\end{equation}

\section{Simulation studies}
In this section we study the empirical performance of our proposed procedure. 
We compare our procedure against
 \texttt{LLR-IK}: 
local linear regressions using the MSE optimal bandwidth~\citep{IK}, \texttt{LLR-RBC}:
local linear regression using the robust bias correction~\citep{RBC}, \texttt{GPR}: Gaussian Process Regressions~\citep{GPR_branson}, \texttt{HRDD}: Local linear regression with a hierarchical model~\citep{sugasawa2023hierarchical}. For \texttt{HRDD} and \texttt{HGPR}, we set a MCMC iteration of 3000, and discard the first 500 iterations as burn-ins.
For the methods which do not accommodate heterogeneity, we apply the procedure to each sub-population individually. For the joint confidence intervals, we combine the individual confidence intervals after adjusting for multiple testing.


We also consider a localized version of \texttt{HGPR}, termed \texttt{HGPR-CUT}, which only considers observations within a window around the threshold $z=0$ instead of all observations. \cite{GPR_branson} also apply local Gaussian process regression to RDD, and note that it improves computation and reduces bias when the stationary covariance assumption is violated. We set our window based on the bandwidth selected by \texttt{LLR-IK}, an approach also taken in \cite{GPR_branson}. More specifically, suppose the bandwidth obtained from \texttt{LLR-IK} is $h_{IK}$, then for \texttt{HGPR-CUT}, we remove any observations that are outside of the interval $[-h_{IK},h_{IK}]$. In the case where the running variable distribution is skewed to the right, we double the window in one direction, $[-h_{IK},2h_{IK}]$, to account for the lack of treated observations, and vice versa.

For each of the data generating procedures below, we perform $500$ replications. For each method we compare the performance of the point estimates by recording the average root mean squared error (RMSE), $\sqrt{\sum_{j=1}^J (\hat{\delta}_{j}-\delta_{j})^2/J}$.
We also compare the performance of confidence/credible intervals by recording average (over all $J$ sub-populations) interval length and empirical coverage.
In Appendix~\ref{app:numerical}, we also report the mean absolute error of the point estimates $\sum_{j=1}^J |\hat{\delta_{j}}-\delta_{j}|$, the simultaneous coverage of joint confidence/credible intervals for $(\delta_j \, :\, j = 1, \ldots, J)$ and the volume of the simultaneous confidence/credible regions.

\subsection{Data generating procedures}
We consider 3 different data generating processes (DGP). In each DGP, we generate $Y_{ij} = f_j(Z_{ij}) + T_{ij}\delta_j +\varepsilon_{ij}$ but vary $f_j$, $\delta_j$, $Z_{ij}$,  and~$\varepsilon_{ij}$. 

DGP1 follows a setting from~\citet{HSU2019468}, where the treatment effect is $0$ for all sub-populations and the conditional mean functions are similar for sub-populations whose indices are close. We set $J= 10$ and $n_j = 100, 200$ for all $j$. 
\begin{equation*}\tag{DGP1}
    \begin{split}
        \delta_j &= 0 \quad \text{ for all } j\\
        f_j(z) &= -0.555-0.0553j + 0.581z + 0.0060jz - 0.058z^2 + 0.01074j^2 \\
        \epsilon_{ij} & \sim N(0,0.1^2)\\
        Z_{ij} & \sim U(-1,1).\\
    \end{split}
\end{equation*}

DGP2 is adapted from \citet{sugasawa2023hierarchical} where the treatment effects are drawn from a Gamma distribution and the conditional mean functions for each sub-population are randomly drawn. Notably, this setting violates our assumption of continuous derivatives at the threshold. 
We set $J = 25, 50$ and $n_j = 100$ for all $j$. 
Due to the computational complexity of \texttt{HGPR}, for DGP2 when $J=50$ we only include \texttt{HGPR-CUT}, not~\texttt{HGPR}. 

\begin{equation*}\tag{DGP2}
    \begin{split}
        f_j(z) &= \begin{cases} a_{j1}z+a_{j2}z^2+a_{j3}z^3, \ z<0 \\
        b_{j1}z+b_{j2}z^2+b_{j3}z^3,\ z\geq 0 \end{cases}\\
        \epsilon_{ij} & \sim N(0,\sigma^2_j)\\
        Z_{ij} & \sim 2\times {\rm Beta}(2,4)-1\\
        \sigma^2_j &\sim U(0.5,1.2) \\
        \delta_j &\sim {\rm Gamma}(3,1)-3 \\
        (a_{j1},b_{j1})  &\sim {\rm Unif}(0.4,1.4), \  a_{j2}  \sim {\rm Unif}(3,7), \ a_{j3}  \sim {\rm Unif}(9,11), \\ b_{j2}  &\sim {\rm Unif}(5,9), \  b_{j3}  \sim {\rm Unif}(3,5). \\
    \end{split}
\end{equation*}

Finally, in DGP3 we intentionally violate many of our model assumptions to examine the robustness of our procedure. We set $J = 10$ and let $n_j = 100, 200$ and sample the data: 
\begin{equation*}\tag{DGP3}
    \begin{split}
        f_j(z) &= a_0+b_1z+a_2z^2+a_3z^3+\sum_{k}c_k(z-a_k)^3\\
        Z_{ij} & \sim U(-1,1)\\
        \sigma^2_j &\sim U(0.25,0.5) \\
        (a_1, a_2, a_3, a_4) & \sim N(0,\sigma_a^2I_4) \\
        c_k& \sim N(0,\sigma_a^2).
    \end{split}
\end{equation*}
Each $f_j$ is a cubic spline with coefficients randomly drawn for each subgroup, and knots $\{a_k\}$ placed at $(-.9, -.8, \ldots, .8, .9)$. Thus, the $f_j$ functions are not infinitely differentiable as would be implied when using the square exponential kernel. To encourage the functions $f_j$ to be highly varied, we set $\sigma_a = 10$. We also consider 2 different distributions for the treatment effects, and 2 different distributions for the random errors:
\begin{itemize}
    \item[(I)]  $\delta \sim N(0,S)$, with $S_{ij} = \rho^{|i-j|}$ for $\rho = .8$,
    \item[(II)] $P(\delta_j =\tau_1) = P(\delta_j =\tau_2) = 1/2$ where $(\tau_1,\tau_2)\sim {\rm Unif}(-3,3)$,
    \item[(A)] $\epsilon_{ij}\sim ({\rm Binom}(5,0.5)-5\times 0.5)/(5 \times 0.25)$,
    \item[(B)] $\epsilon_{ij}\sim \textrm{Gamma}(4,2)-2$.
\end{itemize}
In setting (I), we consider an AR(1) type correlation between treatment effects, where assuming a correlated structure should be beneficial to the hierarchical models {\tt HGPR} and {\tt HRDD}. In setting (II), each sub-population is assigned either treatment effect $\tau_1$ or $\tau_2$. This creates a scenario where shrinkage toward a common mean should be detrimental when $|\tau_1 - \tau_2|$ is large.  Finally, the 2 settings for the random error explores the performance of HGPR when the error is skewed or discrete.

\subsection{Simulation results}

We first note that for DGP1 and DGP3, \texttt{LLR-RBC} crashes in 2\% of the replicates due to an insufficient number of observations within the bandwidth for computing the bias correction term. For DGP 2, this occurs in more than 60\% of replicates, and in the tables below we only report results from the replicates which successfully returned a point estimate and confidence interval. \texttt{HRDD} by default uses the \texttt{LLR-RBC} bandwidth for each subgroup; thus, for DGP2, we use the ``global" setting for \texttt{HRDD}, which fits a bandwidth based on all observations aggregated across subgroups, bypassing the 60\% crash rate of \texttt{LLR-RBC}.

\begin{table}[ht]
\caption{Root Mean Square Error (RMSE)}
\centering
\begin{tabular}{|c|cccccc|}
  \hline
  Setting & HGPR & HGPR-CUT & LLR-IK & LLR-RBC & GPR & HRDD \\ 
  \hline
  DGP1, $n_j=100$ & 0.017 & 0.014 & 0.115 & 0.133 & 0.073 & 0.180  \\ \hline
  DGP1, $n_j=200$ & 0.017 & 0.017 & 0.066 & 0.077 & 0.048 & 0.117  \\ \hline
  DGP2, $J=25$ & 0.330 & 0.280 & 0.715 & 1.786 & 0.594 & 0.450  \\ \hline
  DGP2, $J=50$ &  & 0.372 & 0.704 & 1.794 & 0.588 & 0.418  \\ \hline
  DGP3, (A-I), $n_j=100$ & 0.177 & 0.230 & 1.015 & 0.847 & 0.651 & 0.426  \\ \hline
  DGP3, (A-II), $n_j=100$ & 0.213 & 0.279 & 1.007 & 0.821 & 0.643 & 0.438  \\ \hline
  DGP3, (B-I), $n_j=100$ & 0.175 & 0.227 & 1.059 & 0.852 & 0.665 & 0.439  \\ \hline
  DGP3, (B-II), $n_j=100$ & 0.215 & 0.282 & 1.032 & 0.813 & 0.665 & 0.429  \\ \hline
  DGP3, (A-I), $n_j=200$ & 0.137 & 0.182 & 0.629 & 0.516 & 0.442 & 0.292  \\ \hline
  DGP3, (A-II), $n_j=200$ & 0.148 & 0.205 & 0.623 & 0.516 & 0.437 & 0.306  \\ \hline
  DGP3, (B-I), $n_j=200$ & 0.132 & 0.179 & 0.634 & 0.508 & 0.447 & 0.289  \\ \hline
  DGP3, (B-II), $n_j=200$ & 0.154 & 0.208 & 0.632 & 0.517 & 0.457 & 0.307\\
   \hline
   
\end{tabular}
\label{sim_mse}
\end{table}

Table~\ref{sim_mse} contains the average MSE for each procedure. The RMSE for \texttt{HGPR} is typically smaller than the RMSE of \texttt{HGPR-CUT}; this may be because \texttt{HGPR} utilizes all the data, thus exhibiting less variability, whereas \texttt{HGPR-CUT} only considers a subset of the data. Nonetheless, the RMSE for both \texttt{HGPR} and \texttt{HGPR-CUT} are drastically smaller than the other procedures. Even \texttt{HRDD}, which is specifically tailored for the setting we consider, has an MSE which is at least twice as large as \texttt{HGPR-CUT}.

\begin{table}[ht]
\caption{Average coverage of 95\% credible/confidence intervals}
\centering
\begin{tabular}{|c|cccccc|}
  \hline
   Setting & HGPR & HGPR-CUT & LLR-IK & LLR-RBC & GPR & HRDD \\ 
  \hline
  DGP1, $n_j=100$ & 0.99 & 0.99 & 0.91 & 0.92 & 0.93 & 1.00  \\ \hline
  DGP1, $n_j=200$ & 0.99 & 0.99 & 0.92 & 0.92 & 0.95 & 1.00  \\ \hline
  DGP2, $J=25$ & 0.98 & 0.95 & 0.88 & 0.92 & 0.98 & 1.00  \\ \hline
  DGP2, $J=50$ &  & 0.95 & 0.88 & 0.92 & 0.98 & 0.99  \\ \hline
  DGP3, (A-I), $n_j=100$ & 0.95 & 0.95 & 0.80 & 0.96 & 0.96 & 0.99  \\ \hline
  DGP3, (A-II), $n_j=100$ & 0.94 & 0.95 & 0.81 & 0.97 & 0.96 & 0.99  \\ \hline
  DGP3, (B-I), $n_j=100$ & 0.95 & 0.95 & 0.81 & 0.96 & 0.96 & 0.99  \\ \hline
  DGP3, (B-II), $n_j=100$ & 0.94 & 0.95 & 0.81 & 0.96 & 0.96 & 0.99  \\ \hline
  DGP3, (A-I), $n_j=200$ & 0.92 & 0.94 & 0.83 & 0.95 & 0.96 & 0.99  \\ \hline
  DGP3, (A-II), $n_j=200$ & 0.93 & 0.95 & 0.83 & 0.95 & 0.96 & 1.00  \\ \hline
  DGP3, (B-I), $n_j=200$ & 0.93 & 0.95 & 0.83 & 0.95 & 0.95 & 0.99  \\ \hline
  DGP3, (B-II), $n_j=200$ & 0.92 & 0.95 & 0.83 & 0.94 & 0.95 & 1.00\\
   \hline
   
\end{tabular}
\label{sim_cover}
\end{table}

In Tables~\ref{sim_cover} and~\ref{sim_length} we show the empirical coverage and length of the confidence/credible intervals.
Both \texttt{HGPR} and \texttt{HGPR-CUT} achieve nominal coverage in DGP1 and even DGP2 where the assumption of infinitely many continuous derivatives is violated. In addition, the credible intervals are substantially shorter than the confidence/credible intervals of the competing methods. In the adversarial setting of DGP3, we see that the credible interval lengths for \texttt{HGPR} and \texttt{HGPR-CUT} are still drastically shorter than the competing methods. However, although \texttt{HGPR-CUT} maintains nominal coverage, the empirical coverage of \texttt{HGPR} dips below the nominal level when $n_j = 200$. This may be because the increased sample size makes the model misspecification in the \texttt{HGPR} model more apparent, whereas the local nature of \texttt{HGPR-CUT} suffers less from model misspecification.


\begin{table}[ht]
\caption{Average interval length of 95\% credible/confidence intervals}
\centering
\begin{tabular}{|c|cccccc|}
  \hline
Setting & HGPR & HGPR-CUT & LLR-IK & LLR-RBC & GPR & HRDD \\ 
  \hline
  DGP1, $n_j=100$ & 0.08 & 0.08 & 0.37 & 0.44 & 0.27 & 2.72  \\ \hline
  DGP1, $n_j=200$ & 0.06 & 0.06 & 0.23 & 0.28 & 0.18 & 2.34  \\ \hline
  DGP2, $J=25$ & 1.56 & 1.16 & 2.25 & 5.53 & 2.78 & 2.25  \\ \hline
  DGP2, $J=50$ &  & 1.43 & 2.24 & 5.67 & 2.78 & 2.15  \\ \hline
  DGP3, (A-I), $n_j=100$ & 0.69 & 0.90 & 2.38 & 5.69 & 2.63 & 2.70  \\ \hline
  DGP3, (A-II), $n_j=100$ & 0.81 & 1.05 & 2.37 & 5.59 & 2.64& 2.61  \\ \hline
  DGP3, (B-I), $n_j=100$ & 0.69 & 0.90 & 2.39 & 6.25 & 2.62 & 2.71  \\ \hline
  DGP3, (B-II), $n_j=100$ & 0.80 & 1.04 & 2.38 & 5.56 & 2.63 & 2.60  \\ \hline
  DGP3, (A-I), $n_j=200$ & 0.49 & 0.69 & 1.56 & 2.48 & 1.77 & 2.09  \\ \hline
  DGP3, (A-II), $n_j=200$ & 0.55 & 0.81 & 1.55 & 2.47 & 1.76 & 2.02  \\ \hline
  DGP3, (B-I), $n_j=200$ & 0.48 & 0.69 & 1.56 &2.49 & 1.76 & 2.08  \\ \hline
  DGP3, (B-II), $n_j=200$ & 0.55 & 0.78 & 1.56 &2.48 & 1.78 & 2.02  \\ \hline
   \hline
   
\end{tabular}
\label{sim_length}
\end{table}

\section{Study of Party Advantages in the U.S. Senate}
Political scientists have investigated whether a political party enjoys an advantage in elections where the current incumbent is from that party; see, e.g.,~\citet{gelman1990estimating, levitt1997decomposing, lee2008randomized, erikson2015using, caughey2017elections}. Some studies have also noted that the incumbency advantage has changed over time: at times increasing or decreasing~\citep{cox1996did, jacobson2015personal}. 

We use our method to analyze electoral incumbency advantages in U.S. Senate elections using data from \cite{senate}. 
The data set includes 1201 elections from 1914 to 2010 across all 50 states, and we consider separate time periods as different sub-populations with potentially heterogeneous treatment effects. Each U.S. Senate election is an observation unit, the treatment $T_{ij}$ is whether or not the incumbent is from the Democratic party and the outcome $Y_{ij}$ is the Democratic party's margin of victory (i.e., the Democratic candidate's vote share minus the vote share of the Republican party candidate) 
 in the current election.  
The running variable $Z_{ij}$ is the Democratic party's margin of victory in the previous election  and the cutoff is 0.  Since each U.S. senator serves a 6 year term, we define the previous election as the election 6 years prior in the same state.

\begin{figure}
\centering
        \includegraphics[totalheight=9.5cm]{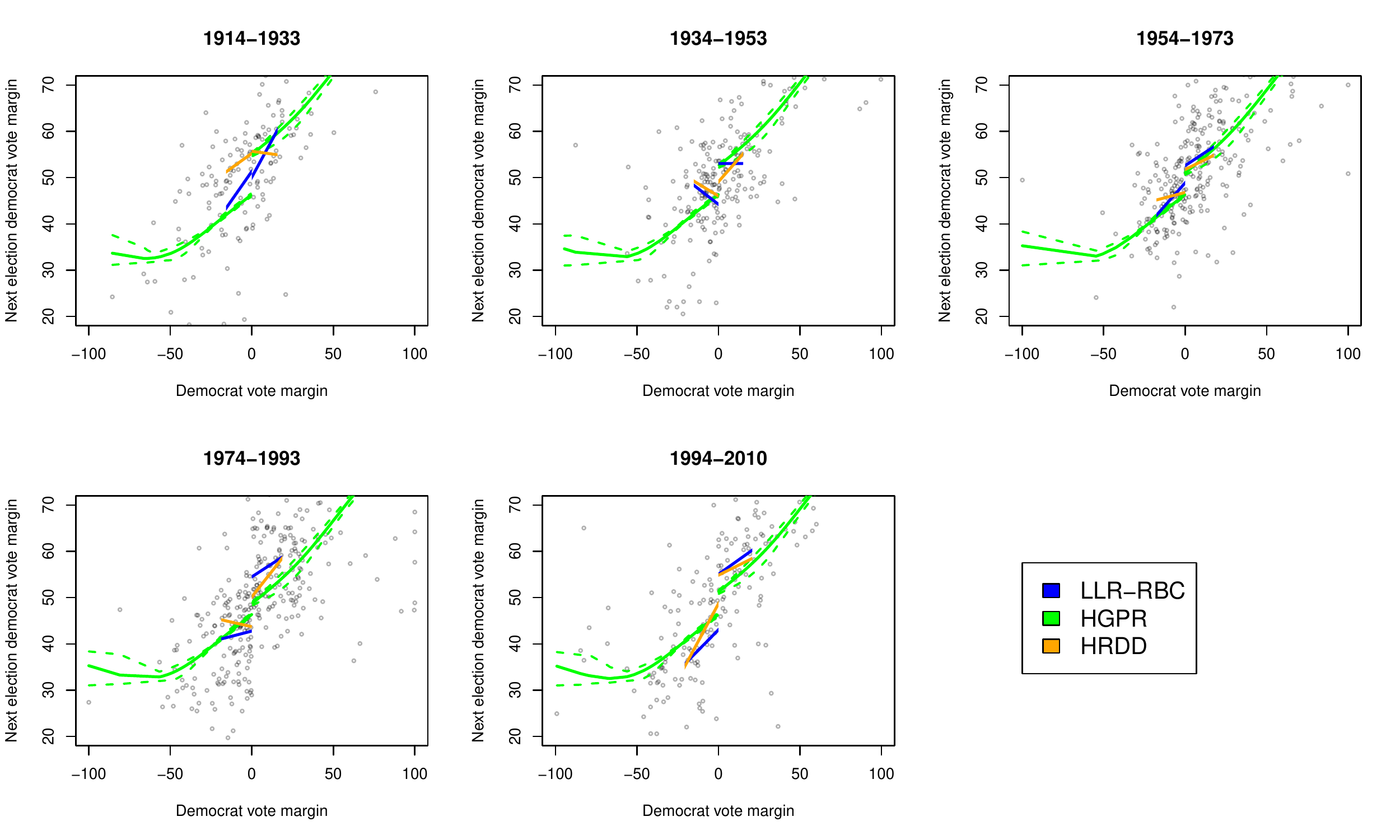}
    \caption{Conditional expectations fit by HGPR and the local linear fit of HRDD and LLR plotted on top of the data points of each group. The solid curve of HGPR is the posterior mean of the conditional mean function, with the dotted curves being the 95\% credible bands. The discontinuous jump of HGPR is the posterior mean of the treatment effect. The HRDD curves are computed using the posterior mean of the local slope and intercept. }
    \label{senate_scatter}
\end{figure}

\citet{senate} use this dataset to examine the incumbent party advantage, and---depending on the model specification---they find a treatment effect around $7.4-9.4\%$. Here we ask the question: Does the incumbency advantage change from 1914 to 2010? To answer this question, we divide the data set into 5 subgroups, each being a roughly 20 year period. We choose 5 because the group-specific sample size ranges between 170 to 312, allowing \texttt{LLR-RBC} to run without problems. We compare \texttt{HGPR}, \texttt{HRDD}, and \texttt{LLR-RBC} when applied to the data set. Both \texttt{HGPR} and \texttt{HRDD} were run with 5000 MCMC iterations, and the first 1000 iterations were discarded as burn-ins. The MCMC chains from both \texttt{HGPR} and \texttt{HRDD} were confirmed to be unimodal and stationary after discarding the burn-ins. Trace plots and density plots from the \texttt{HGPR} chain can be found in Appendix ~\ref{senate:diagnostics}. Table~\ref{senate_table} shows the point estimator and lengths of 95\% marginal confidence intervals for \texttt{HGPR}, \texttt{LLR-RBC}, and \texttt{HRDD}. As expected, the estimates of the treatment effect for  \texttt{HGPR} and \texttt{HRDD} both exhibit shrinkage towards a common mean and result in a smaller range of the point estimates and shorter interval lengths.  
The large confidence interval lengths for \texttt{LLR-RBC} suggests that simply applying robust LLR to each group is ineffective due to inflation of the noise-to-signal ratio when we subdivide the data set. 
To further assess how our estimates differ from \texttt{HRDD} and \texttt{LLR-RBC}, we plotted the fitted conditional mean function, or the local linear fit on top of the scatter plot of the data in Figure~\ref{senate_scatter}. From Figure~\ref{senate_scatter}, we observe that both \texttt{HRDD} and \texttt{LLR-RBC} have a tendency to produce local linear estimates with sharp slopes, which could indicate high sensitivity to the choice of bandwidth. Notably, the \texttt{HRDD} estimate has a negative slope on one side of the conditional mean function for year groups 1914--1933, 1934--1963, and 1984--2003, despite the scatter plot showing a positive trend across all groups. These time periods also happen to be the groups where the \texttt{HGPR} point estimate differs greatly from that of \texttt{HRDD}. \texttt{HGPR} on the other hand captures the positive trend across all groups and is visibly more stable. 

Finally, we apply the credible region method from Eq.~\eqref{post_elip} to test the sharp null hypothesis $H_0: \delta_j=0, j=1,\dots,5$, and the homogeneous null hypothesis $H_0: \delta_1=\dots=\delta_5 = C, \text{for some $C$}$. In both tests, we reject the null hypothesis and thus conclude that the incumbency advantage treatment effect varies across time periods.

\begin{table}[ht]
\caption{The incumbency advantage point estimate and 95\% confidence/credible interval for each time period 
}
\centering
\begin{tabular}{|c|c|c|c|}
    \hline
Time period  & HGPR & HRDD & LLR   \\ \hline 
    1914-1933 & 8.68 (4.39,12.95) & 0.57 (-1.30,2.49) & -2.80 (-10.97,5.36)  \\ \hline
    1934-1953 & 6.19 (2.51,10.07) & 3.01 (0.36,5.84) & 10.2 (2.14,18.15)  \\ \hline
    1954-1973 & 4.69 (1.12,8.39) & 5.02 (3.56,6.32) & 3.30 (-1.47,8.07)  \\ \hline
    1974-1993 & 2.36 (-1.20,6.33) & 6.30 (4.97,7.87) & 11.85 (6.05,17.65)  \\ \hline
    1994-2010 & 4.90 (0.79,9.27) & 5.57 (2.73,8.59) & 11.48 (-1.76,24.72)  \\ 
    \hline
\end{tabular}
\label{senate_table}
\end{table}

\section{Discussion}
In this paper, we propose a procedure for estimating heterogeneous treatment effects using a RDD with known sub-populations. Specifically, we employ a Bayesian hierarchical Gaussian process regression nonparametric regression model which allows for borrowing information across sub-populations.
We can sample from the posterior with Metropolis-Hasting steps within a Gibbs sampler, and---under mild conditions---the posterior concentrates around the true values. Most notably, we show empirically using extensive numerical  simulations that our method  outperforms existing methods. The point estimates have substantially lower mean squared errors, and the credible interval lengths are much shorter while maintaining nominal frequentist coverage across a wide range of settings. 

Interesting future directions include adapting the procedure to the fuzzy RDD setting using the complier and non-complier setup discussed in \citet{chib_greenberg_simoni_2022}. In addition, the model could potentially incorporate pre-treatment covariates by adding additional layers of hierarchy. 


\FloatBarrier

\appendix
\newpage

\section{Proofs}
We now prove Theorem~\ref{thm:hgpr_consistency} by verifying that HGPR satisfies the conditions for posterior consistency required by Theorem 1 of \cite{CHOI20071969} which we restate below. 

\begin{theorem}{(Theorem 1 from \citet{CHOI20071969})} \label{bay_consistency}
 Let $\{Y_i\}_{i=1}^\infty$ be independently distributed with densities $\{\eta_i(\cdot;\theta)\}_{i=1}^\infty$ with respect to a common 
 $\sigma$-finite measure, where the parameter $\theta$ belongs to a abstract measurable space $\Theta$. The densities $\eta_i$ are assumed to be jointly measurable. Let $\theta_0\in \Theta$ be the true value, and let $P_{\theta_0}$ be the corresponding joint distribution of $\{Y_i\}_{i=1}^\infty$. Let $\{U_n\}_{n=1}^\infty$ be a sequence of subsets in $\Theta$, and let $\theta$ have prior $\Pi$. Define:
\begin{equation*}
    \begin{split}
        \Lambda_i(\theta_0,\theta) &= \log \frac{\eta_i(Y_i;\theta_0)}{\eta_i(Y_i;\theta)}, \\
        K_i(\theta_0,\theta) &= E_{\theta_0}(\Lambda_i(\theta_0,\theta)), \\
        V_i(\theta_0,\theta) &= \var (\Lambda_i(\theta_0,\theta)). \\
    \end{split}
\end{equation*}
Suppose the following assumptions hold:
\begin{enumerate}[label=(\Alph*)]
\item Prior positivity of neighborhoods:
There exists a set $B$ with $\Pi(B)>0$ such that
\begin{enumerate}[label=(\arabic*)]
    \item[(A1)] $\sum_{i=1}^{\infty}V_i(\theta_0,\theta)/i^2<\infty,\forall \theta\in B$
    \item[(A2)] For all $\epsilon >0$, $\Pi(B\cap \{\theta:K_i(\theta_0,\theta)<\epsilon,\forall i\})>0$.
  \end{enumerate}
\item Existence of tests: There exist test functions $\{\Psi_n\}_{n=1}^{\infty}$, sets $\{\Theta_n\}_{n=1}^{\infty}$, and constants $C_1,C_2,c_1,c_2>0$ such that
\begin{enumerate}[label=(\arabic*)]
    \item[(B1)] $\sum_{n=1}^{\infty}E_{\theta_0}\Psi_n<\infty$
    \item[(B2)] $\sup_{\theta\in U_n^C\cap\Theta_n}E_{\theta}(1-\Psi_n)\leq C_1e^{-c_1n}$
    \item[(B3)] $\Pi(\Theta_n^C)\leq C_2e^{-c_2n}$.
  \end{enumerate}
\end{enumerate}
Then 
\begin{equation*}
    \Pi(U_n^C|Y_{i}, \ldots, Y_n) \rightarrow 0\quad a.s.\quad [P_{\theta_0}].
\end{equation*}
\end{theorem}

\newpage

\subsection{Proof of Theorem \ref{thm:hgpr_consistency}}

\begin{proof}[Proof of Theorem \ref{thm:hgpr_consistency}]
When the $\delta_j$ are independent and both $g$ and the kernel functions are assumed to be fixed and known, observations from different subgroups are independent. This effectively makes the model (with a slight change of notation):
\begin{equation}
    \begin{split}
        & f_j(z=\cdot)\sim GP(g,K_j), \quad j\in[J] \\
        & Y_{ij} = 1_{(Z_{ij}<0)}[f_j(Z_{ij})+\epsilon_{ij}^-]+1_{(Z_{ij}\geq0)}[f_j(Z_{ij})+\delta_j+\epsilon_{ij}^+]\\
        & \epsilon_{ij}^+\sim N(0,\sigma^2_{+j})\\
        & \epsilon_{ij}^-\sim N(0,\sigma^2_{-j})\\
        &\sigma_{+j} \sim \nu_+ \\
        &\sigma_{-j} \sim \nu_- \\
        & \delta_j \sim P_{\delta} \\
    \end{split}
\end{equation}
Recall that the parameters of interest are $\theta = (f_{1:J},\sigma_{+1:J},\sigma_{-1:J},\delta_{1:J})\in \mathcal{H}^J \times (\mathbb{R}^+)^{2J}\times \mathbb{R}^J$
and we have the sets:
\begin{equation*}
    U_{j\epsilon} = \{(f_j,\sigma_{+j},\sigma_{-j},\delta_j):d_{P_Z}(f_j,f_{j0})<\epsilon,|\sigma_{+j}/\sigma_{+j0}-1|<\epsilon,|\sigma_{-j}/\sigma_{-j0}-1|<\epsilon,|\delta_j-\delta_{j0}|<\epsilon\},
\end{equation*}
and $U_{\epsilon} = \bigtimes\limits_{j=1}^J U_{j\epsilon}$, where $A\bigtimes B$ denotes the Cartesian product of the sets $A$ and $B$. Because $J$ is fixed, if the marginal posterior distribution for each $j=1,\dots,J$ satisfies
\begin{equation}\label{gpr_consistency}
\Pi(U^C_{j\epsilon}|Y_{1:n_j,j},Z_{1:n_j,j})\xrightarrow{a.s.\ P_{\theta_0}}0,
\end{equation}
then
\begin{equation*}
\Pi(U^C_{\epsilon}|\B{Y},\B{Z})\xrightarrow{a.s.\ P_{\theta_0}}0.
\end{equation*}

Therefore, it suffices to show posterior consistency for a single group $j$ under a Gaussian process setup for RDD.
Lemma \ref{GPR_A} and \ref{GPR_B} below show that under Assumptions \ref{asm:kernel} and \ref{asm:prior}, the Gaussian process setup satisfies Assumption A and B in Theorem \ref{bay_consistency}. Thus we have that Eq.~\eqref{gpr_consistency} holds, and the proof is complete.
\end{proof}

\subsection{Lemmas}\label{sec:proof_lemma}
In this section we will work with the following Gaussian process model for estimating the RD treatment effect $\delta$ for a single group.
\begin{equation}\label{gpr}
    \begin{split}
        & f(z=\cdot)\sim GP(g,R) \\
        & Y_{i} = 1_{\{Z_{i}<0\}}[f(Z_{i})+\epsilon_{i}^-]+1_{\{Z_{i}\geq0\}}[f(Z_{i})+\delta+\epsilon_{i}^+]\\
        & \epsilon_{i}^+\sim N(0,\sigma^2_{+})\\
        & \epsilon_{i}^-\sim N(0,\sigma^2_{-})\\
        &\sigma_{+} \sim \nu_+ \\
        &\sigma_{-} \sim \nu_- \\
        & \delta \sim P_{\delta} \\
    \end{split}
\end{equation}
The parameters of interest are $\theta = (f,\sigma_{+},\sigma_{-},\delta)\in \mathcal{H} \times (\mathbb{R}^+)^{2}\times \mathbb{R}$. We will make the following assumptions throughout this section:
\begin{enumerate}
    \item \textbf{GPR Assumption 1}: The running variable $Z_{ij}\stackrel{\rm iid}{\sim} P_Z$, with $P_Z$ being defined on [$-$1,1]
    \item \textbf{GPR Assumption 2}: The Gaussian process $f$ has a continuously differentiable mean function $g$ and the kernel function $R(z,z')$ has continuous 4th partial derivatives. In addition, $\nu_+,\nu_-,P_{\delta}$ assigns positive probability to every neighborhood of the true values $\sigma_{+0},\sigma_{-0},\delta_{0}$, and the distribution $P_{\delta}$ is subexponential; i.e. there exist $K>0$ such that $P(|\delta|\geq t) = O(e^{-Kt})$.
\end{enumerate}
We will assume the function $g$ and the kernel $R$ to be known and fixed. Without loss of generality, we assume $g=0$. The above assumptions are essentially the same as those of Theorem \ref{thm:hgpr_consistency}, except without the group structure. The following Lemmas will show that assumptions A and B from Theorem \ref{bay_consistency} hold for the model in Eq~\eqref{gpr}, thus showing posterior consistency for the GPR RDD model.

\begin{lemma}\label{GPR_A}
Let $\theta_0$ denote the true value and define the set \[B(\gamma) = \left\{(f,\sigma_-,\sigma_+,\delta) \;:\; \left \Vert f-f_0\right \Vert_{\infty}< \gamma,\left \vert\frac{\sigma_-}{\sigma_{-0}} - 1\right \vert < \gamma,\left \vert\frac{\sigma_+}{\sigma_{+0}}-1\right \vert < \gamma,|\delta-\delta_0|<\gamma\right\}.\]
The set $B(1/2)$ satisfies Assumption A from Theorem~\ref{bay_consistency}.
\end{lemma}
\begin{proof}[Proof of Lemma \ref{GPR_A}]
We first introduce the notation $P_+ = P(Z_i\geq 0)$, and $P_- = 1-P_+$. When $\theta$ is fixed, we also use the more compact notation $\Lambda_i = \Lambda_i(\theta_0,\theta)$; we define $K_i$ and $V_i$ analogously. In addition, we let $C$ denote some constant which may depend on the parameters $\theta$ but does not depend on $n$; the specific value of $C$ may also change from line to line.

By GPR Assumption 2, $\Pi(B(\gamma)) > 0$ for any $\gamma  > 0$ as long as \[\mathcal{H}(\left\{f: \left \Vert f-f_0\right \Vert_{\infty}< \gamma \right\})>0,\] 
where $\mathcal{H}(\cdot)$ is the Gaussian process prior. Indeed, this is guaranteed by Theorem 4.2 from \citet{tokdar}; thus, $\Pi(B(1/2)) > 0$.

We now verify (A2). We first note that $B(r) \subset B(s)$ for any $0\leq r < s$. In addition, the distribution of $Y_i$ is $ 1_{\{Z_i\geq 0\}}N(f(Z_i),\sigma^2_+)+1_{\{Z_i<0\}}N(f(Z_i)+\delta,\sigma^2_-)$. Thus, conditional on the sign of $Z_i$ the density of $Y_i$ is a normal density so that
\begin{equation}
    \begin{split}
        K_i &= \E\left[\E_{\theta_0}[\Lambda_i|Z_i]\right] \\
        &= \E[\E_{\theta_0}[\Lambda_i|Z_i,Z_i\geq 0]|Z_i\geq 0]P_+ + \E[\E_{\theta_0}[\Lambda_i|Z_i,Z_i< 0]|Z_i< 0]P_-.
        \end{split}
    \end{equation}
Furthermore, for any $\theta \in B(\gamma)$ for $0 < \gamma < 1/2$ we have for all $i$:
\begin{equation}
    \begin{split}
        \E[\E_{\theta_0}[\Lambda_i|Z_i,Z_i\geq 0]|Z_i\geq 0] &= \frac{1}{2}\log\left(\frac{\sigma^2_+}{\sigma^2_{+0}}\right) - \frac{1}{2}\left(1-\frac{\sigma^2_+}{\sigma^2_{+0}}\right) \\&+ \frac{1}{2}\int_{0}^1\frac{(f_0(z_i)-f(z_i)+\delta_0-\delta)^2}{\sigma^2_+}dP_Z \\
        &\leq \gamma/\sigma^2_{+0}.
    \end{split}
\end{equation}
An analogous result can be shown for when $Z_i<0$. Thus, letting $\gamma(\epsilon) = \min\{1/2, \sigma^2_{+0}\epsilon\}$, we have
\[B(\gamma(\epsilon)) \subseteq \left\{\theta : K_i(\theta_0, \theta) < \epsilon \; \forall i\right\}. \]
This implies that for any $\epsilon >0$
\begin{equation}
\begin{aligned}
        \Pi\left(B(1/2)\cap \{\theta:K_i(\theta_0,\theta)<\epsilon,\forall i\}\right) & \geq \Pi\left(B(1/2)\cap B(\gamma(\epsilon))\right)\\
    &= \Pi\left(B(\gamma(\epsilon))\right) > 0,
\end{aligned} 
\end{equation}
which proves (A2)

To show (A1), we use the law of total variance. 
\begin{equation}
    \begin{split}
        V_i &= \var_{\theta_0}(\Lambda_i) \\
        &= \var_{\theta_0}(\Lambda_i|Z_i<0)P_- + \var_{\theta_0}(\Lambda_i|Z_i\geq0)P_+ \\
        &- 2P_-P_+\E_{\theta_0}(\Lambda_i|Z_i<0)\E_{\theta_0}(\Lambda_i|Z_i\geq0) \label{eq:Vi} \\
        &+ \E_{\theta_0}(\Lambda_i|Z_i<0)^2P_-P_+ + \E_{\theta_0}(\Lambda_i|Z_i\geq0)^2P_-P_+ .
    \end{split}
\end{equation}
Note that $\E_{\theta_0}(\Lambda_i|Z_i<0) \leq \gamma/\sigma^2_{-0}$ and $\E_{\theta_0}(\Lambda_i|Z_i\geq0) \leq \gamma/\sigma^2_{+0}$. For the terms in the second row, note that we can apply the law of total variance again:
\begin{equation}
    \begin{split}
        \var_{\theta_0}(\Lambda_i|Z_i\geq0) &= \E[\var_{\theta_0}(\Lambda_i|Z_i,Z_i\geq0)|Z_i\geq0]+ \var_{\theta_0}[\E_{\theta}(\Lambda_i|Z_i,Z_i\geq0)|Z_i\geq0] \\
        = 2(-\frac{1}{2}&+\frac{\sigma^2_{+0}}{2\sigma^2_+})^2+\int_{0}^1 \left[\frac{\sigma^2_{+0}}{\sigma^2_+}(f(z_i)-f_0(z_i)+\delta-\delta_0)\right]^2dP_Z < 5\gamma^2. 
    \end{split}
\end{equation}
The second equality comes from Section 4.2.1 of \citep{CHOI20071969}. A similar argument can be applied to $Z_i<0$, allowing us to conclude that $V_i$ is uniformly bounded and $\sum_i V_i/i^2 < \infty$. Thus $B(1/2)$ satisfies assumption A in Theorem \ref{bay_consistency}.

\end{proof}

\begin{lemma}\label{GPR_B}
Fix some $0< \epsilon<1$ and define  
\[
    U_{\epsilon} = \{(f,\sigma_{+},\sigma_{-},\delta):d_{P_Z}(f,f_{0})<\epsilon,|\sigma_{+}/\sigma_{+0}-1|<\epsilon,|\sigma_{-}/\sigma_{-0}-1|<\epsilon,|\delta-\delta_{0}|<\epsilon\},\]
    and let $U_n=U_{\epsilon}$ for all $n$. Then there exists test functions $\{\Phi_n\}_{n=1}^{\infty}$, and sets $\{\Theta_n\}_{n=1}^{\infty}$, such that assumption B from Theorem~\ref{bay_consistency} holds.
\end{lemma}
\begin{proof}[Proof of Lemma \ref{GPR_B}]
Let $M_n =n^{3/4}$ and define the set $\Theta_n = \Theta_{1n}\times \Theta_{2n} \times (\mathbb{R}^+)^2$, where
\begin{equation}
    \begin{split}
        \Theta_{1n} &= \{f: \Vert f\Vert_{\infty}<M_n, \Vert f'\Vert _{\infty}<M_n\}\\
        \Theta_{2n} &= [-n,n].
    \end{split}
\end{equation}

Let $\tilde \Theta_{1n}(t)$ denote a $t$-net of $\Theta_{1n}$ with respect to the $\Vert \cdot \Vert_\infty$ norm and let
$\tilde \Theta_{2n}(t)$ denote a $t$-net of $\Theta_{2n}$. Furthermore, let $N_{1t}$ and $N_{2t}$ denote the cover numbers of $\Theta_{1n}$ and $\Theta_{2n}$ respectively. Then by \citet[Theorem 2.7.1]{vaart-wellner}, we have for fixed $t >0$ that $\log N_{1t}=O(M_n)$. In addition, we have $N_{2t} =O(n)$. 

We will construct a test for each element of $\tilde \Theta_{1n}(t) \times \tilde \Theta_{2n}(t)$ and ultimately combine all the tests for a single test. Specifically, let $t = \min\{\epsilon/2,r/4,s/4\}$, with $0<r<\min(\epsilon^2/4,4\sigma_{-0}\sqrt{\epsilon-\epsilon^2})$ and $0<s<\min(\epsilon/2,$ $4\sigma_{+0}\sqrt{\epsilon-\epsilon^2})$. 

We first verify (B3). With our definition of $\Theta_n$, we have $\Pi(\Theta_n^C) \leq \mathcal{H}(\Theta_{1n}^C)+P_{\delta}(\Theta_{2n}^C)$.
For the choice of $M_n$, theorem 5 of \cite{ghosal_roy} gives $\mathcal{H}(\Theta_{1n}^C)\leq C_1 \exp(-C_2n)$. By properties of a subexponential distribution, we also have $P_{\delta}(\Theta_{2n}^C)\leq C_3\exp(-C_4n)$. Thus we have $\Pi(\Theta_n^C) = O(e^{-cn})$.

To verify the remainder of Assumption B, we will use $n_+$ and $n_-$ to denote the number of treated and control observations. Since $n_+$ follows a binomial distribution with probability $P_+$, we can use Hoeffding's inequality to obtain $P(n_+\leq nP_+/2)\leq \exp(-n(P_+/2)^2)$. Using the same argument for $n_- > P_{-}/2$, we can conclude that  $\min\{n_+, n_-\} > C n$ with probability at least $1 - 2\exp(-n\min\{P_+^2, P_-^2\})$. Thus, in the remainder of this section, we will assume that there exists a positive constant $C$ such that $n_- > C n$ and $n_+ > C n$ for all $n$.

Let $(f^{l}, \delta^{l'})$ denote some element of $\tilde \Theta_{1n}(t) \times \tilde \Theta_{2n}(t)$. We define $f_{li} = f^{l}(z_i)$ and $f_{0i} = f_{0}(z_i)$ for $i=1,\dots,n$.  Let $\gamma>0$, $c_n = n^{3/7}$, $b_i = 1_{\{f_{li}\geq f_{0i}\}} - 1_{\{f_{li}< f_{0i}\}}$, and $b = 1_{\{\delta^{l'}\geq \delta_{0}\}} - 1_{\{\delta^{l'}< \delta_{0}\}}$. We further define the following indicator tests:
\begin{equation}
    \begin{split}
        \Psi_{1n}[f^{l},\delta^{l'},\gamma] &= 1\left\{\sum_{i:z_i<0}b_i\left(\frac{Y_j-f_{0i}}{\sigma_{-0}}\right)>2c_n\sqrt{n_-}\right\} \\
        \Psi_{2n}[f^{l},\delta^{l'},\gamma] &= 1\left\{\sum_{i:z_i<0}\bigg(\frac{Y_j-f_{0i}}{\sigma_{-0}} \bigg)^2 >n_-(1+\gamma) \right\}\\
        \Psi_{3n}[f^{l},\delta^{l'},\gamma] &= 1\left\{\sum_{i:z_i<0}\bigg(\frac{Y_j-f_{li}}{\sigma_{-0}} \bigg)^2 <n_-(1-\gamma^2) \right\}\\
        \Psi_{4n}[f^{l},\delta^{l'},\gamma] &= 1\left\{\sum_{i:z_i\geq0}b\bigg(\frac{Y_j-f_{0i}-\delta_0}{\sigma_{+0}}\bigg)>2c_n\sqrt{n_+}\right\} \\
        \Psi_{5n}[f^{l},\delta^{l'},\gamma] &= 1\left\{\sum_{i:z_i\geq0}\bigg(\frac{Y_j-f_{0i}-\delta_0}{\sigma_{+0}}\bigg)^2>n_+(1+\gamma) \right\} \\
        \Psi_{6n}[f^{l},\delta^{l'},\gamma] &= 1\left\{\sum_{i:z_i\geq0}\bigg(\frac{Y_j-f_{0i}-\delta^{l'}}{\sigma_{+0}}\bigg)^2<n_+(1-\gamma^2) \right\}.\
        \end{split}
\end{equation}
  
We will define the test:
\begin{equation}
    \Psi_n[f^{l},\delta^{l'},\gamma] = \max_{k\in[6]}\Psi_{kn}[f^{l},\delta^{l'},\gamma].
\end{equation}

We now verify (B1) by analyzing the Type I error of the proposed test functions. We fix $\gamma>0$ and we observe that $\E_{P_0}\Psi_n[f^{l},\delta^{l'},\gamma]\leq \sum_{i=1}^6\E_{P_0}\Psi_{in}[f^{l},\delta^{l'},\gamma]$. There exists a positive constant $D$ such that $n_- > D n$ and $n_+ > Dn$ for all $n$. Therefore, Theorem 2 from \citep{CHOI20071969} guarantees $\E_{P_0}\Psi_{1n}[f^{l},\delta^{l'},\gamma]$, $\E_{P_0}\Psi_{2n}[f^{l},\delta^{l'},\gamma]$, and $\E_{P_0}\Psi_{3n}[f^{l},\delta^{l'},\gamma]$ are all order $O(e^{-2c_n^2})$. Since we have assumed the running variables to be fixed, any $Y_j$ that are in the summations of  $\Psi_{4n}[f^{l},\delta^{l'},\gamma],\ \Psi_{5n}[f^{l},\delta^{l'},\gamma], \ \Psi_{6n}[f^{l},\delta^{l'},\gamma]$ necessarily follow the distribution $N(\delta_0+f_0(z_j),\sigma^2_{+0})$ under the null hypothesis. Furthermore, replacing $b_i$ with $b$ does not change the standard normal distributions involved in $\Psi_{4n}$, therefore, we can apply Theorem 2 again to obtain $\E_{P_0}\Psi_{4n}[f^{l},\delta^{l'},\gamma],$ $\ \E_{P_0}\Psi_{5n}[f^{l},\delta^{l'},\gamma], \ \E_{P_0}\Psi_{6n}[f^{l},\delta^{l'},\gamma] = O(e^{-2c_n^2})$, which together imply $\E_{P_0}\Psi_n[f^{l},\delta^{l'},\gamma]= O(e^{-2c_n^2})$. 
We will define the test function $\Psi_n = \max_{f^{l}\in \tilde \Theta_{1n}(t),\delta^{l'}\in \tilde \Theta_{1n}(t)} \Psi_n[f_j,\delta^{l'},\epsilon/2]$. Because $\log(N_{1t}N_{2t}) = o(c_n^2)$, this results in the Type I error:
\begin{equation}
    \begin{split}
        \E_{p_0}\Psi_n &\leq \sum_{f^{l}\in \tilde \Theta_{1n}(t),\delta^{l'}\in \tilde \Theta_{1n}(t)}\E_{p_0}\Psi_n[f^{l},\delta^{l'},\epsilon/2] \\
        &\leq CN_{1t}N_{2t}e^{-2c_n^2} = O(e^{-c_n^2}).
    \end{split}
\end{equation}
Thus we have verified (B1). 

We now verify (B2). Let $f$, $\sigma^2_{+}$, $\sigma^2_{-0}$, and $\delta$ denote the true values. We first note that 
$\E_P(1-\Psi_n[f^{l},\delta^{l'},\epsilon/2])\leq \min_{k\in [6]}\E_P(1-\Psi_{kn}[f^{l},\delta^{l'},\epsilon/2])$, with $P$ denoting the joint distribution of $\{Y_i\}_{i=1}^{n}$ under $\theta = (f,\sigma_+,\sigma_-,\delta)\in U_{\epsilon}^C\cap \Theta_n$. Furthermore, the Type II error of $\Psi_n$ is no larger than the minimum of the individual Type II error of each $\Psi_n[f^{l},\delta^{l'},\epsilon/2]$ test. Thus, we only need to find one pair $(f^{l},\delta^{l'})$ with exponentially small Type II error. 

Due to our choice of $t$ for the $t$-net, we know that for any $f \in \Theta_{1n}$, there exists $f^{l}$ such that $\Vert f^{l}-f\Vert_{\infty}<t$; and if $d_{P_Z}(f,f_0)>\epsilon$, we have $d_{P_Z}(f^{l},f_0)>\epsilon/2$, and according to Lemma 11 in \citep{CHOI20071969}:
\begin{equation}\label{lem:choi11}
    P\Big(\sum_{i:z_i<0}|f^{l}(Z_j)-f_0(Z_j)|\geq rn_-\Big) > 1-e^{-Cn}, \ r\in (0,\epsilon^2/4).
\end{equation}
This shows that our choice of $r$ satisfies the event in Eq~\eqref{lem:choi11} with high probability. For $s$, we go through a similar argument. For any $\delta \in \Theta_{2n}$, there exists $\delta^{l'}$ such that $|\delta^{l'}-\delta|<t$; and if $|\delta_0-\delta|>\epsilon$, we have $|\delta^{l'}-\delta_0|>\epsilon/2 > s$. We will let this pair of $(f^{l},\delta^{l'})$ be the ones used for the derivation of Type II error. Also, note that with this specific pair, we have $\Vert f-f^{l}\Vert_{\infty}<r/4$ and $|\delta-\delta^{l'}|<s/4$. We now verify that the Type II error is exponentially small on this pair.

In the following analysis, we will condition on the running variable $Z$, as well as both $\sum_{i:z_i<0}^n|f_{li}-f_{0i}|>rn_-$ and $|\delta^{l'}-\delta_0|>s$. We will show that the conditional type II error is exponential small, and appeal to Eq.~\eqref{lem:choi11} to obtain the unconditional type II error bounds. We consider the null hypothesis:
\begin{equation}
    H_0: f=f_0,\ \sigma_-=\sigma_{-0},\ \sigma_+=\sigma_{+0},\ \delta=\delta_0.
\end{equation}
There is a total of 15 different possible alternatives, yet whenever $d_{P_Z}(f,f_0)>\epsilon$ or $|\sigma_-/\sigma_{-0}-1|>\epsilon$ occurs in the alternative, we can use Theorem 2 from \citep{CHOI20071969} to obtain exponential bounds on the Type II error by operating on $E_P(1-\Psi_{1n}[\epsilon/2]),\ E_P(1-\Psi_{2n}[\epsilon/2])$ and $E_P(1-\Psi_{3n}[\epsilon/2])$. Therefore, we are only left with 3 possible alternatives:
\begin{itemize}
    \item $H_A: f=f_0,\sigma_-=\sigma_{-0},|\sigma_+/\sigma_{+0}-1|>\epsilon,|\delta-\delta_0|>\epsilon$.
    \item $H_A: f=f_0,\sigma_-=\sigma_{-0},|\sigma_+/\sigma_{+0}-1|>\epsilon,\delta=\delta_0$.
    \item $H_A: f=f_0,\sigma_-=\sigma_{-0},\sigma_+=\sigma_{+0},|\delta-\delta_0|>\epsilon$.
\end{itemize}
We can rewrite these 3 alternatives into the following cases:
\begin{enumerate}
    \item $|\delta-\delta_0|>\epsilon,\sigma_+\leq(1+\epsilon)\sigma_{+0}$.
    \item $\sigma_+>(1+\epsilon)\sigma_{+0}$.
    \item $\sigma_+<(1-\epsilon)\sigma_{+0}$.
\end{enumerate}
All 3 cases assume $f=f_0,\sigma_-=\sigma_{-0}$. For the case $|\delta-\delta_0|>\epsilon,\sigma_+\leq(1+\epsilon)\sigma_{+0}$, we will assume $n_+$ to be large enough so that $c_n/\sqrt{n_+}<s/(4\sigma_{+0})$, then we have:
\begin{equation}
    \begin{split}
        E_P(1-\Psi_n[f^{l},\delta^{l'},\epsilon/2]) &\leq E_P(1-\Psi_{4n}[\epsilon/2])\\
        &= P\left(\sum_{i:z_i\geq0}b(Y_j-f_{0i}-\delta_0)/\sigma_{+0}\leq 2c_n\sqrt{n_+}\right)\\
        &= P\biggl\{\frac{1}{\sqrt{n_+}}\sum_{i:z_i\geq0}b\left(\frac{Y_j-f_{0i}-\delta}{\sigma_+}\right)+\frac{1}{\sqrt{n_+}}\sum_{i:z_i\geq0}b\left(\frac{\delta-\delta^{l'}}{\sigma_+}\right) \\
        &+\frac{1}{\sqrt{n_+}}\sum_{i:z_i\geq0}\left|\frac{\delta^{l'}-\delta_0}{\sigma_+}\right|\leq 2c_n\frac{\sigma_{+0}}{\sigma_+}\biggl\} \\
        &\leq P\left\{\frac{1}{\sqrt{n_+}}\sum_{i:z_i\geq0}b\left(\frac{Y_j-f_{0i}-\delta}{\sigma_+}\right)\leq \frac{s\sqrt{n_+}}{4\sigma_+}- \frac{s\sqrt{n_+}}{\sigma_+} +2c_n\frac{\sigma_{+0}}{\sigma_+}  \right\} \\
        &\leq \Phi\left(-\frac{s\sqrt{n_+}}{4\sigma_{+0}(1+\epsilon)}\right) =O\left(\frac{1}{\sqrt{n}}e^{-Cn}\right) = O\left(e^{-Cn}\right).
    \end{split}
\end{equation}
The last line is due to Mill's Inequality for the standard normal distribution.
For the case $\sigma_{+}>(1+\epsilon)\sigma_{+0}$, we denote $W\sim \chi^2_{n_+}$ and $W'$ be a noncentral $\chi^2_{n_+}$ with noncentrality parameter $n_+(\delta-\delta_0)^2/\sigma_+^2$, then we use Chernoff bound as follows:
\begin{equation}
    \begin{split}
        E_P(1-\Psi_n[f^{l},\delta^{l'},\epsilon/2]) &\leq E_P(1-\Psi_{5n}[\epsilon/2])\\
        &\leq P\left( \sum_{i:z_i\geq0}\left( \frac{Y_j-f_{0i}-\delta_0}{\sigma_+} \right)^2\leq\frac{\sigma^2_{+0}}{\sigma^2_+}n_+(1+\epsilon) \right)\\
        &= P\left( W'\leq \frac{\sigma^2_{+0}}{\sigma^2_+}n_+(1+\epsilon)\right) \\
        &\leq  P\left( W\leq \frac{\sigma^2_{+0}}{\sigma^2_+}n_+(1+\epsilon)\right) \\
        &\leq P\left( W\leq n_+/(1+\epsilon)\right) \\
        &\leq \exp(-n_+t/(1+\epsilon))(1-2t)^{-n_+/2}, \ \forall t<0 \\
        &=\exp\left( -n_+\frac{\epsilon^2-\epsilon^3}{4(1+\epsilon)} \right), \quad \text{by letting } t=-\epsilon/2.
    \end{split}
\end{equation}
Finally, case 3 assumes $\sigma_+<(1-\epsilon)\sigma_{+0}$. We will make use of $\Psi_{6n}$:
\begin{equation}
    \begin{split}
        E_P(1-\Psi_n[f^{l},\delta^{l'},\epsilon/2]) &\leq E_P(1-\Psi_{6n}[\epsilon/2])\\
       & = P\left( \sum_{i:z_i\geq0}\left( \frac{Y_j-f_{0i}-\delta^{l'}}{\sigma_{+0}} \right)^2\geq n_+(1-\epsilon^2/4)\right)\\
       &= P\left( W'\geq n_+(1-\epsilon^2/4)\frac{\sigma^2_{+0}}{\sigma^2_+}\right),
    \end{split}
\end{equation}
where $W'$ is a noncentral $\chi^2_{n_+}$, with noncentrality parameter $d = \frac{n_+(\delta - \delta^{l'})^2}{\sigma_+^2}$. $W'$ has the moment generating function $M_{W'}(t) = (1-2t)^{-n_+/2}\exp\{-d/2[1-(1-2t)^{-1}]\}$ for all $t<1/2$. Therefore, we can use Chernoff bound to show:
\begin{equation}
    \begin{split}
        &P\left( W'\geq n_+(1-\epsilon^2/4)\frac{\sigma^2_{+0}}{\sigma^2_+}\right)\\
        &\leq \exp\left\{ \frac{n_+}{2}\left [-\log(1-2t)-\left( 1-\frac{1}{1-2t}\right)\frac{d}{n}-2t(1-\epsilon^2/4)\frac{\sigma^2_{+0}}{\sigma^2_+} \right ]\right\}\\
        &\leq \exp\left\{ \frac{\sigma^2_{+0}tn_+}{\sigma^2_+}\left [\frac{1}{1-2t}\left( (1-\epsilon)^2+\frac{s^2}{16\sigma^2_{+0}}\right)-(1-\epsilon^2/4) \right ]\right\}, \ \text{since }|\delta^{l'}-\delta|<s/4\\
        &\leq \exp\left\{ \frac{\sigma^2_{+0}tn_+}{\sigma^2_+}\left [\frac{1}{1-2t} (1-\epsilon)-(1-\epsilon^2/4) \right ]\right\}, \ \text{since }s^2<16\sigma^2_{+0}(\epsilon-\epsilon^2)\\
        &\leq \exp\left\{ -n_+t^*\frac{\sigma^2_{+0}\epsilon^3}{4\sigma^2_+}\right\}, \ \text{by setting }\frac{1}{1-2t^*} = \frac{1-(1+\epsilon)\epsilon^2/4}{1-\epsilon}\\
        &\leq \exp\left\{ -n_+t^*\frac{\epsilon^3}{4(1-\epsilon)^2}\right\}
    \end{split}
\end{equation}
Since we have shown that the test with $(f^{l},\delta^{l'})$ has exponentially small Type II error $\Psi_n$ also has exponentially small Type II error. Thus (B2) is verified.
\end{proof}

\newpage

\subsection{Identifiability}
\begin{proof}[Proof of Lemma \ref{lem:id}]
    Our identifiability conditions are sub-population-specific versions of the conditions in \cite{identify}. We will maintain the same definition for $T_j^+,T_j^-$ as in Lemma \ref{lem:id}. We will similarly define:
    \begin{equation}
        \begin{split}
            Y_j^- &= \lim_{z\rightarrow 0^-}\E(Y_{ij}|Z_{ij}=z) \\
            Y_j^+ &= \lim_{z\rightarrow 0^+}\E(Y_{ij}|Z_{ij}=z).
        \end{split}
    \end{equation}
    
Under the assumption of Lemma \ref{lem:id}, and recallig that $Y_{ij} = Y_{ij}(0)+\delta_{ij}T_{ij}$, we can observe that:
\begin{equation}
    \begin{split}
        Y^+_j-Y^-_j &= \lim_{z\rightarrow 0^+}\E(Y_{ij}|Z_{ij}=z) - \lim_{z\rightarrow 0^-}\E(Y_{ij}|Z_{ij}=z) \\
        &= \lim_{z\rightarrow 0^+}[\E(T_{ij}\delta_{ij}|Z_{ij}=z)-\E(T_{ij}\delta_{ij}|Z_{ij}=-z)]\\
        &+\lim_{z\rightarrow 0^+}[\E(Y_{ij}(0)|Z_{ij} = z)-\E(Y_{ij}(0)|Z_{ij} = -z)]\\
        &= \E(\delta_{ij}|Z_{ij} = 0)[T^+_j-T^-_j].\\
    \end{split}
\end{equation}
Thus, $\delta_j = \frac{Y^+_j-Y^-_j}{T^+_j-T^-_j}$ for all $j=1,\dots,J$.
\end{proof}

\newpage
\section{Additional Numerical Results}
\label{app:numerical}
To the measure performance of point estimates, we also includes measurement of the absolute bias averaged across subgroups, as well as mean absolute error (MAE), calculated as: $\sum_{j=1}^J |\hat{\delta}_{j}-\delta_{j}|/J$. To measure inferential performance, we include the multivariate coverage (Multi-Cover), which is calculated using Bonferroni correction for competing methods, and using the posterior ellipsoid method for HGPR and HGPR-CUT. Finally, we also record the $J^{th}$ root of the volume for the confidence region ($V^{1/J}$). From the numerical results, we observe that both versions of HGPR obtain a significantly smaller bias and MAE compared to any existing methods. In terms of coverage, although the Multi-coverage of HGPR and HGPR-CUT is occasionally slightly less than the targeted 95\%, the volume of their confidence regions are drastically smaller than those of existing methods.

\begin{table}[ht]
\caption{Simulation result for DGP1 and 2}
\centering
\begin{tabular}{|c|c|ccccccc|}
  \hline
   Settings & Method & Length & Cover & $|$Bias$|$ & RMSE & MAE & Multi-Cover & $V^{1/J}$ \\ 
  \hline
  \multirow{5}{*}{\shortstack{DGP1\\ J=10\\ $n_j=100$}}&
   GPR & 0.267 & 0.929 & 0.042 & 0.073 & 0.058 & 0.867 & 0.371 \\ &
   LLR IK & 0.372 & 0.913 & 0.004 & 0.115 & 0.080 & 0.802 & 0.483 \\&
   LLR RBC & 0.444 & 0.916 & 0.004 & 0.133 & 0.093 & 0.821 & 0.577 \\&
   HRDD & 2.718 & 1.000 & 0.039 & 0.180 & 0.141 & 1.000 & 3.801\\&
   HGPR & 0.078 & 0.993 & 0.006 & 0.017 & 0.010 & 0.988 & 0.107 \\&
   HGPR-CUT & 0.078 & 0.993 & 0.004 & 0.014 & 0.010 & 0.998 & 0.111\\
  \hline
  \multirow{5}{*}{\shortstack{DGP1\\ J=10\\ $n_j=200$}}&
   GPR & 0.184 & 0.946 & 0.025 & 0.048 & 0.038 & 0.911 & 0.256 \\ &
   LLR-IK & 0.232 & 0.917 & 0.002 & 0.066 & 0.051 & 0.816 & 0.321 \\ &
   LLR-RBC & 0.275 & 0.920 & 0.002 & 0.077 & 0.059 & 0.827 & 0.381 \\ &
   HRDD & 2.335 & 1.000 & 0.041 & 0.117 & 0.093 & 1.000 & 3.270 \\&
   HGPR & 0.060 & 0.995 & 0.005 & 0.017 & 0.008 & 0.978 & 0.083 \\ &
   HGPR-CUT & 0.057 & 0.994 & 0.003 & 0.017 & 0.008 & 0.998 & 0.081 \\
  \hline
  \multirow{3}{*}{\shortstack{DGP2\\ J=25}}&
   GPR & 2.778 & 0.976 & 0.148 & 0.594 & 0.466 & 0.963 & 4.213\\ &
   LLR-IK & 2.245 & 0.879 & 0.062 & 0.715 & 0.554 & 0.550 & 3.399 \\ &
   HRDD & 2.251 & 0.988 & 0.012 & 0.450 & 0.352 & 1.000 & 3.383\\&
    HGPR & 1.562 & 0.982 & 0.008 & 0.330 & 0.261 & 0.988 & 1.845 \\&
   HGPR-CUT & 1.115 & 0.952 & 0.016 & 0.280 & 0.221 & 0.923 & 1.241\\
   \hline
   \multirow{3}{*}{\shortstack{DGP2\\ J=50}}&
   GPR & 2.783 & 0.978 & 0.148 & 0.588 & 0.462 & 0.953 & 4.335\\ &
   LLR-IK & 2.243 & 0.880 & 0.060 & 0.704 & 0.544 & 0.402 & 3.613 \\ &
   HRDD & 2.147 & 0.990 & 0.019 & 0.418 & 0.324 & 0.990 & 3.381 \\&
   HGPR-CUT & 1.434 & 0.947 & 0.011 & 0.372 & 0.289 & 0.920 & 1.744\\
   \hline
   
\end{tabular}
\label{DGP12}
\end{table}

\begin{table}[ht]
\caption{Simulation result for HGPR on DGP3}
\centering
\begin{tabular}{|c|ccccccc|}
  \hline
 Settings & Length & Cover & $|$Bias$|$ & RMSE & MAE & Multi-Cover & $V^{1/J}$\\  
  \hline
(A-I), $n_j=100$ & 0.692 & 0.945 & 0.004 & 0.177 & 0.139 & 0.946 & 0.539 \\ 
(A-II), $n_j=100$ & 0.813 & 0.936 & 0.006 & 0.213 & 0.165 & 0.944 & 0.562 \\ 
(B-I), $n_j=100$ & 0.685 & 0.948 & 0.003 & 0.176 & 0.140 & 0.956 & 0.530 \\ 
(B-II), $n_j=100$ & 0.804 & 0.939 & 0.010 & 0.215 & 0.167 & 0.952 & 0.564 \\ 
(A-I), $n_j=200$ & 0.485 & 0.921 & 0.011 & 0.137 & 0.110 & 0.946 & 0.376 \\ 
(A-II), $n_j=200$ & 0.548 & 0.926 & 0.005 & 0.148 & 0.117 & 0.952 & 0.402 \\ 
(B-I), $n_j=200$ & 0.484 & 0.931 & 0.002 & 0.131 & 0.105 & 0.956 & 0.371 \\ 
(B-II), $n_j=200$ & 0.550 & 0.920 & 0.007 & 0.154 & 0.122 & 0.938 & 0.393 \\ 
   \hline
\end{tabular}
\end{table}

\begin{table}[ht]
\caption{Simulation result for HGPR-CUT on DGP3}
\centering
\begin{tabular}{|c|ccccccc|}
  \hline
 Settings & Length & Cover & $|$Bias$|$ & RMSE & MAE & Multi-Cover & $V^{1/J}$\\  
  \hline
(A-I), $n_j=100$ & 0.901 & 0.953 & 0.008 & 0.230 & 0.180 & 0.954 & 0.756 \\ 
(A-II), $n_j=100$ & 1.053 & 0.953 & 0.007 & 0.279 & 0.214 & 0.946 & 0.822 \\ 
(B-I), $n_j=100$ & 0.891 & 0.953 & 0.005 & 0.227 & 0.179 & 0.956 & 0.761 \\ 
(B-II), $n_j=100$ & 1.036 & 0.948 & 0.004 & 0.282 & 0.217 & 0.964 & 0.808 \\ 
(A-I), $n_j=200$ & 0.687 & 0.940 & 0.017 & 0.182 & 0.143 & 0.938 & 0.583 \\ 
(A-II), $n_j=200$ & 0.808 & 0.949 & 0.012 & 0.206 & 0.162 & 0.948 & 0.633 \\ 
(B-I), $n_j=200$ & 0.690 & 0.952 & 0.004 & 0.179 & 0.141 & 0.956 & 0.588 \\ 
(B-II), $n_j=200$ & 0.778 & 0.948 & 0.006 & 0.208 & 0.159 & 0.965 & 0.613 \\ 
   \hline
\end{tabular}
\end{table}

\begin{table}[ht]
\caption{Simulation result for GPR on DGP3}
\centering
\begin{tabular}{|c|ccccccc|}
  \hline
 Settings & Length & Cover & $|$Bias$|$ & RMSE & MAE & Multi-Cover & $V^{1/J}$\\  
  \hline
(A-I), $n_j=100$  & 2.628 & 0.957 & 0.025 & 0.651 & 0.509 & 0.946 & 3.698 \\ 
(A-II), $n_j=100$ & 2.637 & 0.961 & 0.018 & 0.644 & 0.504 & 0.959 & 3.711 \\ 
(B-I), $n_j=100$ & 2.618 & 0.957 & 0.028 & 0.664 & 0.511 & 0.939 & 3.671 \\ 
(B-II), $n_j=100$  & 2.626 & 0.957 & 0.020 & 0.665 & 0.505 & 0.945 & 3.686 \\ 
(A-I), $n_j=200$ & 1.766 & 0.959 & 0.013 & 0.441 & 0.350 & 0.970 & 2.498 \\ 
(A-II), $n_j=200$ & 1.756 & 0.957 & 0.014 & 0.438 & 0.348 & 0.946 & 2.482 \\ 
(B-I), $n_j=200$ & 1.761 & 0.947 & 0.012 & 0.448 & 0.350 & 0.937 & 2.487 \\ 
(B-II), $n_j=200$ & 1.777 & 0.946 & 0.015 & 0.456 & 0.353 & 0.931 & 2.491 \\
   \hline
\end{tabular}
\end{table}

\begin{table}[ht]
\caption{Simulation result for LLR-IK on DGP3}
\centering
\begin{tabular}{|c|ccccccc|}
  \hline
 Settings & Length & Cover & $|$Bias$|$ & RMSE & MAE & Multi-Cover & $V^{1/J}$\\  
  \hline
(A-I), $n_j=100$ & 2.380 & 0.802 & 0.037 & 1.016 & 0.732 & 0.540 & 3.296 \\ 
(A-II), $n_j=100$ & 2.373 & 0.810 & 0.039 & 1.004 & 0.727 & 0.534 & 3.292 \\ 
(B-I), $n_j=100$ & 2.392 & 0.807 & 0.031 & 1.056 & 0.749 & 0.524 & 3.331 \\ 
(B-II), $n_j=100$ & 2.382 & 0.807 & 0.041 & 1.031 & 0.738 & 0.474 & 3.317 \\ 
(A-I), $n_j=200$ & 1.556 & 0.830 & 0.017 & 0.628 & 0.466 & 0.620 & 2.184 \\ 
(A-II), $n_j=200$& 1.553 & 0.827 & 0.022 & 0.623 & 0.461 & 0.598 & 2.182 \\ 
(B-I), $n_j=200$ &  1.564 & 0.831 & 0.020 & 0.633 & 0.466 & 0.638 & 2.201 \\ 
(B-II), $n_j=200$ &  1.564 & 0.826 & 0.035 & 0.632 & 0.467 & 0.626 & 2.202 \\
   \hline
\end{tabular}
\end{table}

\begin{table}[ht]
\caption{Simulation result for LLR-RBC on DGP3}
\centering
\begin{tabular}{|c|ccccccc|}
  \hline
 Settings & Length & Cover & $|$Bias$|$ & RMSE & MAE & Multi-Cover & $V^{1/J}$\\  
  \hline
(A-I), $n_j=100$ & 5.688 & 0.959 & 0.029 & 0.847 & 0.595 & 0.884 & 6.983 \\ 
(A-II), $n_j=100$ & 5.593 & 0.968 & 0.032 & 0.821 & 0.584 & 0.919 & 6.883 \\ 
(B-I), $n_j=100$ & 6.251 & 0.958 & 0.205 & 0.851 & 0.602 & 0.893 & 6.981 \\ 
(B-II), $n_j=100$ & 5.558 & 0.960 & 0.026 & 0.813 & 0.594 & 0.903 & 6.907 \\ 
(A-I), $n_j=200$ & 2.484 & 0.947 & 0.020 & 0.516 & 0.395 & 0.904 & 3.282 \\ 
(A-II), $n_j=200$ & 2.471 & 0.951 & 0.017 & 0.516 & 0.396 & 0.898 & 3.272 \\ 
(B-I), $n_j=200$ & 2.488 & 0.945 & 0.016 & 0.507 & 0.393 & 0.862 & 3.297 \\ 
(B-II), $n_j=200$ & 2.479 & 0.942 & 0.014 & 0.515 & 0.404 & 0.864 & 3.292 \\ 
   \hline
\end{tabular}
\end{table}

\begin{table}[ht]
\caption{Simulation result for HRDD on DGP3}
\centering
\begin{tabular}{|c|ccccccc|}
  \hline
 Settings & Length & Cover & $|$Bias$|$ & RMSE & MAE & Multi-Cover & $V^{1/J}$\\  
  \hline
(A-I), $n_j=100$ & 2.700 & 0.998 & 0.016 & 0.426 & 0.339 & 1.000 & 3.78 \\ 
(A-II), $n_j=100$ & 2.612 & 0.999 & 0.012 & 0.438 & 0.348 & 1.000 & 3.66 \\ 
(B-I), $n_j=100$ & 2.718 & 0.999 & 0.020 & 0.439 & 0.350 & 1.000 & 3.79 \\ 
(B-II), $n_j=100$ & 2.598 & 0.998 & 0.010 & 0.429 & 0.338 & 1.000 & 3.64 \\ 
(A-I), $n_j=200$ & 2.091 & 1.000 & 0.014 & 0.292 & 0.233 & 1.000 & 2.93 \\ 
(A-II), $n_j=200$ & 2.027 & 1.000 & 0.004 & 0.306 & 0.243 & 1.000 & 2.83 \\ 
(B-I), $n_j=200$ & 2.085 & 1.000 & 0.010 & 0.289 & 0.230 & 1.000 & 2.92 \\ 
(B-II), $n_j=200$ & 2.027 & 0.999 & 0.018 & 0.307 & 0.245 & 0.998 & 2.83 \\ 
   \hline
\end{tabular}
\end{table}

\clearpage
\section{HGPR Convergence Diagnostics}
\label{senate:diagnostics}
Here we present the MCMC trace and density plots for the treatment effects of the five groups in the Senate incumbency application. We also present the plots for the prior standard deviation parameter for the treatments, which is the standard deviation term in the kernel function $K_{\delta}$ in Eq~\eqref{eq:HGPR}. We observe that all parameters are stationary and unimodal, thus suggesting convergence in the MCMC chain.

\begin{figure}[!h]
\centering
\includegraphics[scale=0.95,page=1]{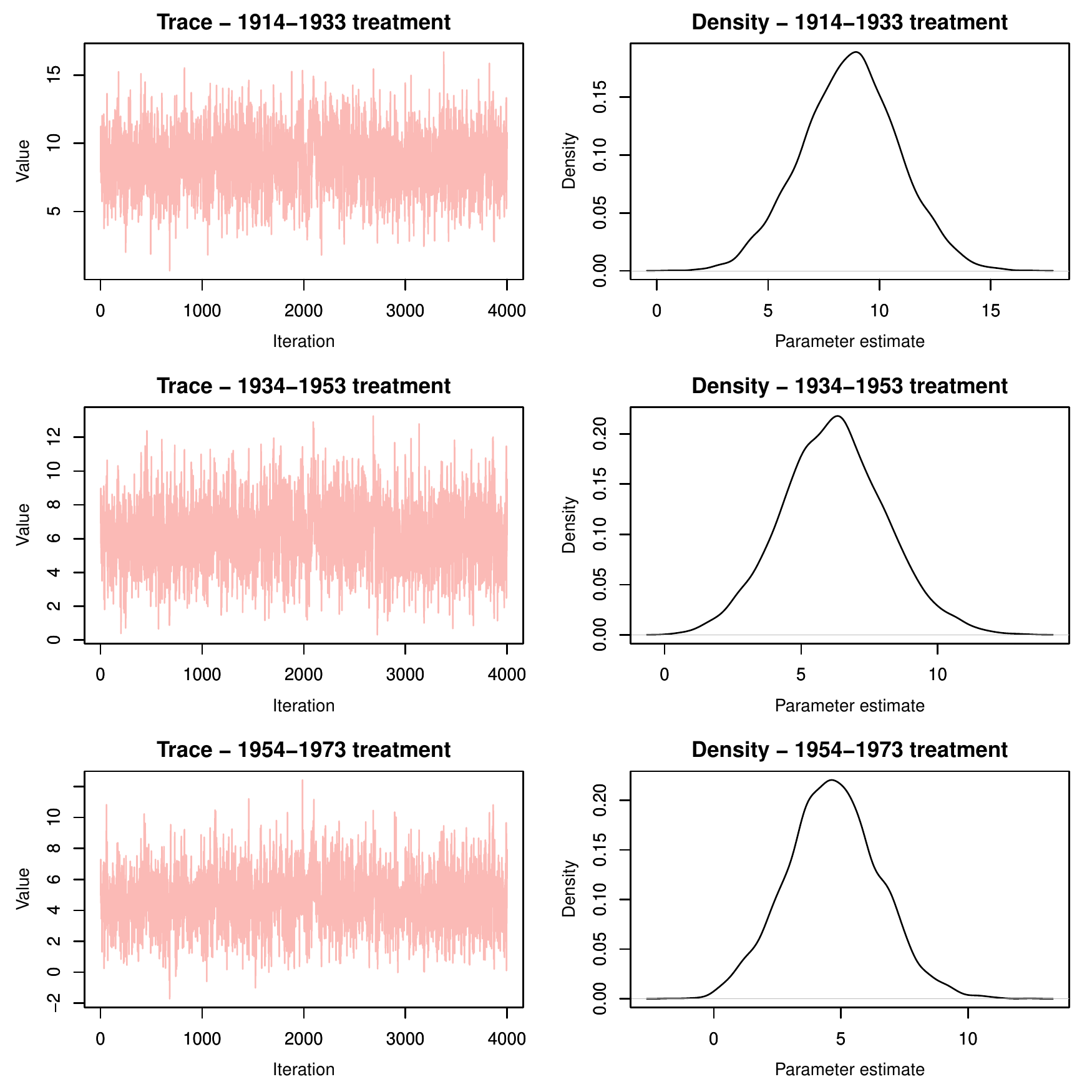}
\label{senate_trace}
\end{figure}

\begin{figure}[!h]
\includegraphics[scale=0.95,page=2]{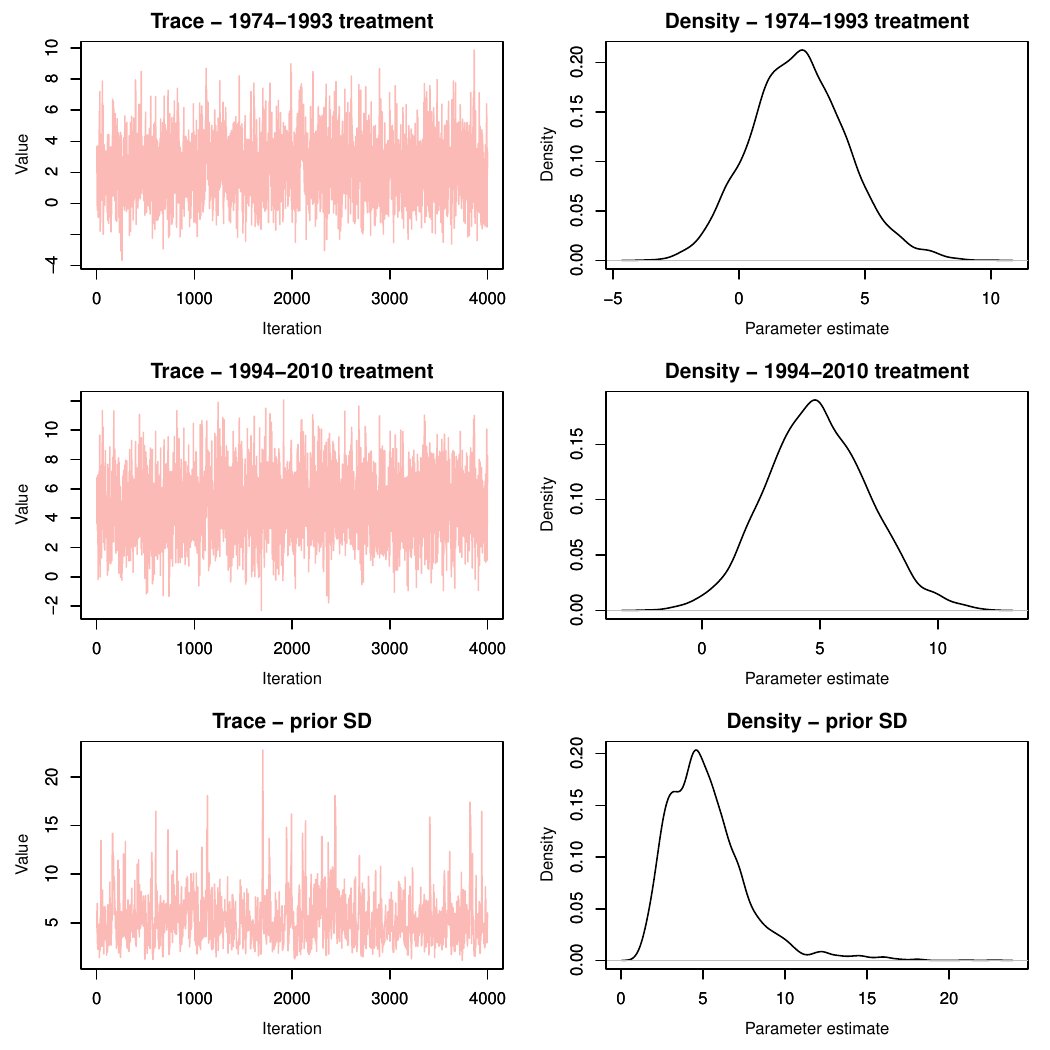}
\end{figure}


\end{document}